\documentclass{article}
\usepackage{amsfonts,amssymb,mathrsfs,graphics,natbib}
\usepackage{latexsym}
\usepackage{amssymb}
\usepackage{amsmath}
\usepackage[english]{babel}
\usepackage{fontenc}
\usepackage{times}
\usepackage[all]{xy}
\usepackage{epsfig}
\usepackage{amsthm}
\usepackage[latin1]{inputenc}
\usepackage{setspace}
\usepackage{yfonts}
\usepackage{lscape}
\usepackage{natbib}
\usepackage{url}
\usepackage{fixltx2e}
\usepackage{alltt}
\usepackage{graphicx}
\usepackage{mathrsfs}
\usepackage{wrapfig}

\graphicspath{{immagini//}}

\usepackage[all]{xy}
\usepackage{setspace}
\usepackage{lscape}

\newcommand{\C}{\mathbb{C}}
\newcommand{\CP}{\mathbb{CP}}

\newcommand{\R}{\mathbb{R}}
\newcommand{\Gr}{\mathbb{G}\mathrm{r}}

\renewcommand{\theenumii}{\arabic{enumii}}
\def\s#1.{\stepcounter{enumii}\theenumi.\theenumii\ \textsl{#1.}\\}

\newcommand\mb{\medbreak}
\newcommand\n{\noindent}
\newcommand\ip{\raise1pt\hbox{\large$\lrcorner$}\,}

\newcommand{\g}{\mathfrak{g}}

\newcommand{\X}{\mathcal{X}}
\newcommand{\Y}{\mathcal{Y}}

\newcommand{\ba}{\begin{array}}\newcommand{\ea}{\end{array}}
\newcommand{\rf}[1]{(\ref{#1})}\def\lie#1({\mathfrak{#1}(}

\newcommand{\nit}[1]{\mb\n\textit{#1.}}

\renewcommand{\ge}{\geqslant}
\renewcommand{\le}{\leqslant}
\newcommand{\la}{\lambda}

\newcommand{\sF}{\mathscr{F}}
\newcommand{\sG}{\mathscr{G}}
\newcommand{\sO}{\mathscr{O}}
\newcommand{\sP}{\mathscr{P}}
\newcommand{\sD}{\mathscr{D}}
\newcommand{\fF}{\mathfrak{F}}
\newcommand{\mV}{\mathcal{V}}
\newcommand{\mH}{\mathcal{H}}

\newcommand{\+}{\!+\!}

\newcommand{\ti}{\!\times\!}
\newcommand{\rr}[1]{[\![#1]\!]}

\def\,{\kern2pt}
\def\.{\,\cdot\,}

\theoremstyle{plain}
\newtheorem{teor}{Theorem}
\newtheorem{prop}[teor]{Proposition}
\newtheorem{lemma}[teor]{Lemma}
\newtheorem{corol}[teor]{Corollary}
\theoremstyle{definition}
\newtheorem{defi}[teor]{Definition}
\theoremstyle{rem}
\newtheorem{rem}[teor]{Remark}

\begin{document}

\title{\bf{Toric moment mappings and Riemannian structures\footnote{With acknowledgement to Geometriae Dedicata. The final publication Toric moment mappings and Riemannian structures DOI: 10.1007/s10711-012-9720-6 is available at www.springerlink.com/content/yn86k22mv18p8ku2/ }}}
\bigskip\author{\bf{Georgi Mihaylov}}

\maketitle

\begin{abstract}\noindent Coadjoint orbits for the group
  $SO(6)$ parametrize Riemannian $G$-reductions in six dimensions, and
  we use this correspondence to interpret symplectic fibrations
  between these orbits, and to analyse moment polytopes associated to
  the standard Hamiltonian torus action on the coadjoint orbits. The
  theory is then applied to describe so-called intrinsic torsion
  varieties of Riemannian structures on the Iwasawa manifold.
\end{abstract}

\section*{Introduction}

In an article of Abbena, Garbiero and Salamon \cite{AGSAHG} (based on
\cite{AGSHGI}), the authors exploit a solid tetrahedron to describe
the set of orthogonal almost complex structures on the Iwasawa
manifold and other 6-dimensional nilmanifolds. Their employment of the
tetrahedron was justified on purely combinatorial grounds, enabling
the 16 Gray--Hervella classes of almost Hermitian structures to be
represented in terms of unions of vertices, faces, edges and other
segments. Our aim is to insert this theory into a more universal
setting.

In the present article, we consider general reductions of a
Riemannian structure specified by a subgroup $G$ of $SO(N)$
stabilizing a $2$-form. The relevance of 2-forms to the intrinsic
torsion of Riemannian structures arises from the isomorphism
 of the space of 2-forms with
the orthogonal Lie algebra ($\Lambda^2T^*_pM
\cong\mathfrak{so}(N)$). Moreover, the orthogonal complement
$\mathfrak{g}^\perp$ of the Lie algebra of $G$ in $\mathfrak{so}(N)$
is a model of the vertical space of the bundle $P/G$ parametrizing
the $G$-structures.

Our point of view emphasizes the natural role that symplectic
geometry plays in the classification of Riemannian structures. The
flag varieties parametrizing reductions of the type we are
considering are merely adjoint (equivalently, coadjoint) orbits of
$SO(6)$. The coadjoint orbits of any compact Lie group $G$ are
precisely the manifolds on which $G$ acts transitively as a compact
group of symplectic automorphisms (see \cite{VGMMCI,GLSSFM}). This
fact enables us to prove that every $G$-structure defined by a
$2$-form on a 6-manifold is associated to a moment mapping from a
flag variety to $\mathfrak{so}(6)^*$. Restriction to a maximum torus
gives rise to a \emph{moment polytope} in $\R^3$, and the shape of
the polytope is sensitive to the exact 2-form chosen. This
construction provides a precise geometrical interpretation of the
tetrahedron introduced in \cite{AGSAHG}.

In our set-up, phenomena involving symplectic fibrations of
coadjoint orbits, such as symplectic quotients and other operations
analysed in \cite{GLSSFM}, have a direct and detailed interpretation
in terms of compatibility conditions for specific Riemannian
structures in six real dimensions. A discussion of certain invariant
subsets in Section~3 enables us to identify the image of torus
moment mappings, and characterize the resulting faces using almost
complex structures. In Section~4, we interpret the well-known Klein
correspondence of projective geometry between elements in $\CP^3$
and $\Gr_2(\C^4)$ from this viewpoint. In the majority of cases,
this and similar correspondences can be clearly visualized in terms
of the moment polytopes characterizing the structures involved.

One of the $G$-structures determined by a $2$-form is the
\emph{mixed structure}, which we define in this article. The theory
that we develop enables one to describe a mixed structure from both
algebraic and geometric point of view in terms of orthogonal almost
complex and almost product structures which are well known.

In fact in the last section we illustrate some implications of the
theory that go beyond mere algebraic and combinatorial aspects.
Since nilmanifolds are parallelizable in a natural way, they provide
a rich source of examples of structures defined globally in terms of
invariant tensors. We describe an application of the theory
involving the classes of various types of Riemannian structures on
the Iwasawa manifold characterized by specific constraints on their
intrinsic torsion. In this sense Corollaries \ref{cor1}, \ref{cor2}
and \ref{cor3} are among the main results of the paper.

\section{SO(6) coadjoint orbits}

A fundamental result guarantees that any orbit of the adjoint action
of a compact Lie group $G$ on its Lie algebra $\g$ intersects the
closure of each Weyl chamber in a single point (see for example
\cite{RBGRLG,ABEM}). This property implies that the set of adjoint
orbits can be parametrized by the closed fundamental Weyl chamber.
The standard identification $\g\cong\g^*$, realized by the Killing
form, also allows us to identify the orbits of the adjoint and the
coadjoint actions.

It is always possible to define a symplectic structure on a
coadjoint orbit $\sO$ such that the inclusion
$\sO\hookrightarrow\g^*$ is the moment map associated to the
Hamiltonian action of $G$ (see \cite{GuSSTP}). This is the
Konstant--Kirillov--Souriau (KKS) structure, defined by
\[ \omega_{\lambda}(\X,\Y)=(\lambda,[X,Y]),\qquad X,Y\in
\mathfrak{g},\]\noindent where $\X,\Y\in T_{\lambda}\sO$ are
determined by the vector fields generated by $X,Y$ i.e.
$\X=\mathrm{ad}_{\lambda}X$ and $\Y=\mathrm{ad}_{\lambda}Y$.

Restricting the group action to the maximal torus $T\subset G$, we
obtain a Hamiltonian torus action on the orbit. The moment map
$\mu_T$ associated to this action consists of the orthogonal
projection to the subalgebra $\mathfrak{t}\subset\g$.  Any coadjoint
orbit $\sO$ intersects $\mathfrak{t}$ in a single orbit of the Weyl
group:
$$\sO\cap\mathfrak{t}=W\cdot\lambda,$$ for some $\lambda\in
\mathfrak{t}$. The points in the intersection of the orbit and
$\mathfrak{t}$ are exactly the points fixed by the action of $T$,
and none of these are found in the interior of the convex polytope
determined by the Weyl orbit of $\lambda$. The celebrated Atiyah and
Guillemin--Sternberg (AGS) Convexity Theorem implies that the image
by $\mu_T$ of an orbit passing through $\lambda\in\mathfrak{t}$ is
the convex hull of the Weyl group orbit of $\lambda$: $$\mu_T
(G\cdot\lambda) = \mathrm{conv}(W\cdot\lambda).$$ See
\cite{Ati,GLSSFM} for more details.

We first apply this general theory to provide a complete description
of the set of $SO(6)$ adjoint orbits. Consider the maximum torus
$T\subset SO(6)$ containing the matrices:

\begin{equation}\label{torus}\left(
                               \begin{array}{ccc}
                                 A_1 & 0 & 0 \\
                                 0 & A_2 & 0 \\
                                 0 & 0 & A_3 \\
                               \end{array}
                             \right),\,\,\,\,A_i=\left(
                                               \begin{array}{cc}
                                                 \cos 2\pi\theta_i&\,\,\, -\sin 2\pi\theta_i \\
                                                  \sin 2\pi\theta_i &\,\,\, \cos 2\pi\theta_i\\
                                               \end{array}
                                             \right),\,\,\,i=1,2,3\end{equation}

The corresponding Lie algebra $\mathfrak{t}\subset \mathfrak{so}(6)$
is generated by:
\begin{equation}\label{torus}\left(
                               \begin{array}{ccc}
                                 B_1 & 0 & 0 \\
                                 0 & B_2 & 0 \\
                                 0 & 0 & B_3 \\
                               \end{array}
                             \right),\,\,\,\,B_i=\left(
                                               \begin{array}{cc}
                                                 0&-\theta_i \\
                                                 \theta_i&0 \\
                                               \end{array}
                                             \right),\,\,\,i=1,2,3\end{equation}

\begin{figure}[!h]
\centering\includegraphics[width=0.6\textwidth]{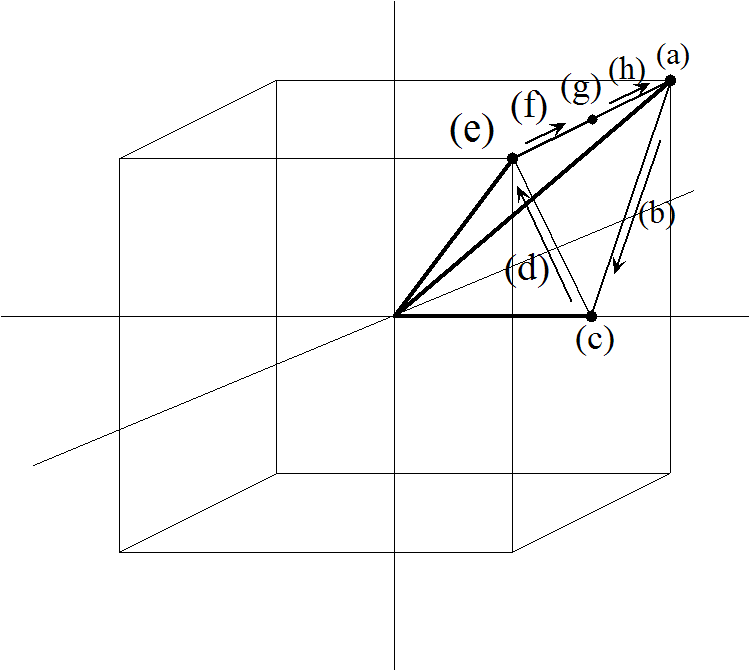}\\
\caption{The fundamental Weyl chamber of $SO(6)$.}
\label{parameter2}\end{figure}

\noindent Using as a basis the matrices $v_i$ such that
$\theta_i=1$, $\theta_j=0$ for $i\neq j$ we can identify
$\mathfrak{t}$ isometrically with $\R^3$. Relative to this basis,
the fundamental weights are:
$$(\theta_1-\theta_2)\, ,\,(\theta_2-\theta_3)\, ,\,(\theta_2+\theta_3),$$
\noindent and the fundamental Weyl chamber $B$ is determined by the
inequalities (see Figure~\ref{parameter2}):
\begin{equation} \theta_1>\theta_2,\quad \theta_2>\theta_3,\quad \theta_2>-\theta_3.
\label{weylsu4}\end{equation} \noindent The bounding cube on the
figure helps to visualize the remaining 23 Weyl chambers in
analogous positions. The Weyl group is generated by the reflections
with respect to the walls of the Weyl chambers.

A \emph{generic} adjoint (equivalently coadjoint) orbit, passing
through an interior point of $B$, is the real 12-dimensional
manifold
$$\sO^{SO(6)}=\frac{SO(6)}{U(1)\times U(1)\times U(1)}\cong\frac{SO(6)}T$$ of ``full'' complex flags of $\R^6$.
A point belonging to the faces or the edges of the closed
fundamental Weyl chamber $\bar{B}$ admits a stabilizer which
contains properly the above torus subgroup, and so gives rise to a
\emph{degenerate} orbit. Points in $\bar{B}$ with the corresponding
stabilizer and orbit are listed in Table 1. We denote by $U(p)$ and
$\tilde{U}(p)$ the unitary subgroups of $SO(2p)$ associated
respectively to the complex structures $J$ and $\tilde{J}$ on
$\R^{2p}$ which act on a standard basis as
 follows:
\[\begin{array}{l}
Je_1=e_2,...,Je_{2p-3}=e_{2p-2},Je_{2p-1}=e_{2p},\\[3pt]
\tilde{J}e_1=e_2,...,\tilde{J}e_{2p-3}=e_{2p-2},\tilde{J}e_{2p-1}=-e_{2p}.
\end{array}\]

\bigbreak\phantom.\hspace{-15pt}

\begin{tabular}{|l|l|l|l|l|l|l|}\hline
                    point in $\bar{B}$       & stabilizer               & orbit                    & image by $\mu_T$                                   &Fig.~\ref{parameter1} \\ \hline\hline
                  $(\alpha,\alpha,\alpha)$   & $U(3)$                   & $\sP^+ \cong \CP^{3}$    & tetrahedron $\Delta_{\sP^+}$                       &(a)\\ \hline\hline
                  $(\alpha,\alpha,-\alpha)$  & $\tilde{U}(3)$           & $\sP^- \cong \CP^{3}$    & tetrahedron $\Delta_{\sP^-}$                       &(e)\\ \hline\hline
                  $(\alpha,0,0)$             & $U(1)\ti SO(4)$          & $\sG\!\cong\!\Gr_2(\R^6)$& octahedron $\Delta_{\sG}$                          &(c)\\ \hline\hline
                  $(\alpha,\beta,\beta)$     & $U(1)\ti U(2)$           & $\sF^+ $                 & {\small truncated tetrahedron} $\Delta_{\sF^+}$    &(b)\\ \hline
                  $(\alpha,\beta,-\beta)$    & $U(1)\ti \tilde{U}(2)$   & $\sF^-$                  & {\small truncated tetrahedron} $\Delta_{\sF^-}$    &(d)\\ \hline
                  $(\alpha,\alpha,\beta)$    & $U(2)\ti U(1)$           & $\sD^+ $                 & {\small skew-cuboctahedron} $\Delta_{\sD^+}$       &(h)\\ \hline
                  $(\alpha,\alpha,0)$        & $U(2)\ti SO(2)$          & $\sD^0 $                 & {\small cuboctahedron} $\Delta_{\sD^0}$            &(g)\\ \hline
                  $(\alpha,\alpha,-\beta)$   & $U(2)\ti \tilde{U}(1)$   & $\sD^-$                  & {\small skew-cuboctahedron} $\Delta_{\sD^-}$       &(f)\\ \hline
\end{tabular}\smallbreak
Table 1 \bigbreak

The general theory (see \cite{ABEM}) leads us to distinguish three
types of orbits. The ``$+$'' orbits are the complex flag manifolds of
$\R^6$ where the total complex structure induces the canonical
orientation on $\R^6$. The ``$-$'' orbits are complex flag manifolds
with a total complex structure inducing the opposite orientation on
$\R^6$. A ``$+$'' orbit is conjugate to a ``$-$'' one by an element in
$O(6)$ exchanging the complex structure.  Then we have the $0$ type
orbits $\sG$ and $\sD^0$ which are partial complex flag manifolds.
From this point of view, the $10$-dimensional orbits reflect
$SO(N)$-inequivalent cases. However in the case under consideration,
since $U(1)$, $\tilde{U}(1)$ and $SO(2)$ are exactly the same group,
the $\sD$-type orbits are all identical.\smallbreak

\begin{figure}[!h]

 \centering\includegraphics[width=0.7\textwidth]{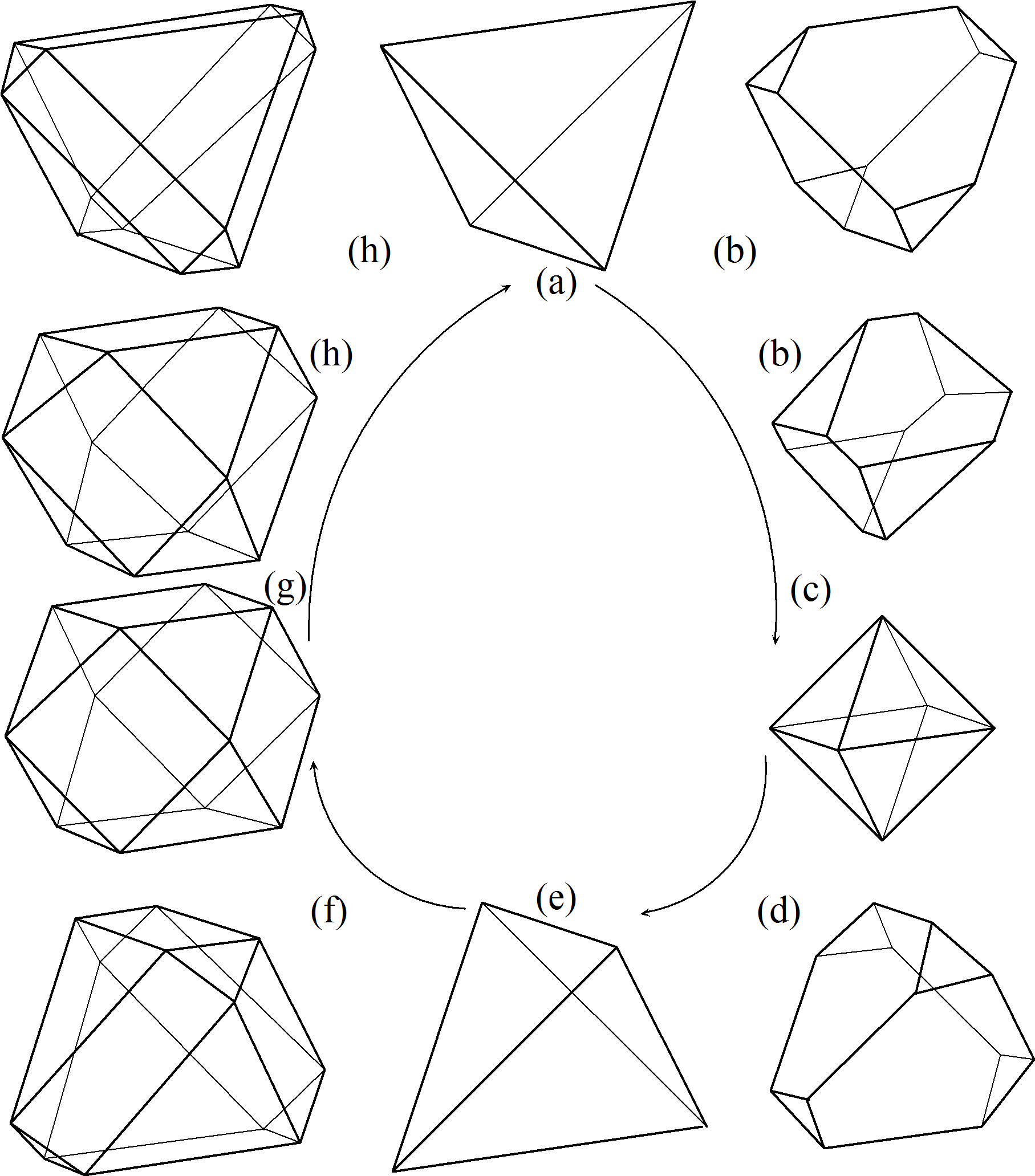}\\
 \caption{Continuous family of moment polytopes for degenerate $SO(6)$
   orbits}
 \label{parameter1}\end{figure}

\begin{figure}[!h]
\centering\includegraphics[width=1.0\textwidth]{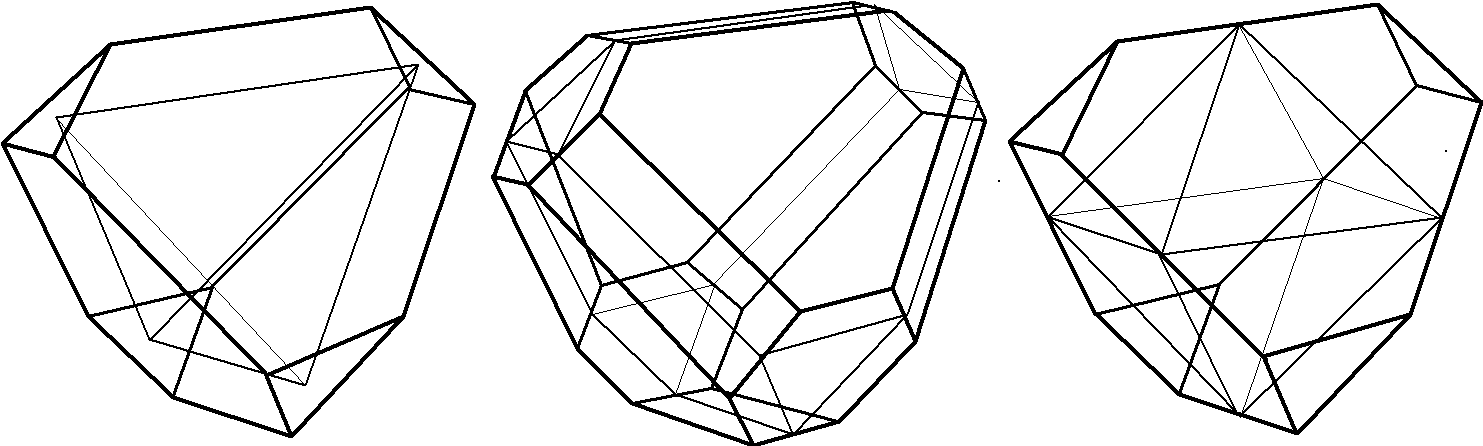}\\
\caption{$\sF^+$ fibres over $\sP^+$; \ $\sO^{SO(6)}$ fibres over
over $\sF^+$; \ $\sF^+$ fibres over $\sG$.}
\label{2fibrazioni}\end{figure}

The Weyl group of $SO(6)$ acts by permuting the coordinates and/or
changing an even number of signs, so the Weyl orbit of
$(\alpha,\alpha,\alpha)$ in $\R^3\cong\mathfrak{t}^*$ consists of
itself and
$(\alpha,-\alpha,-\alpha)$, $(-\alpha,\alpha,-\alpha)$ ,$(-\alpha,-\alpha,\alpha)$.
The resulting tetrahedron is denoted by $\Delta_{\sP^+}$. Analogously
the Weyl orbit of $(\alpha,\alpha,-\alpha)$ is given by itself and the
points $(-\alpha,\alpha,\alpha)$, $(\alpha,-\alpha,\alpha)$ and
$(-\alpha,-\alpha,-\alpha)$, leading to $\Delta_{\sP^-}$. The specific
position of a point in $\bar{B}$ (Figure~\ref{parameter2}) determines
in an analogous way the degeneracy of the Weyl orbit and hence the
polytopes illustrated in Figure~\ref{parameter1}.  Following the
arrows 1--4 we see how the position inside $\bar{B}$ affects the
precise shape of the moment polytope. The image $\mu_T(\sO^{SO(6)})$
is determined by the non-degenerate Weyl orbit represented by the
dodecahedron $\Delta_\sO$ in the middle of Figure~\ref{2fibrazioni}.

\medbreak Obvious inclusion relations between the stabilizer groups
listed in Table 1 show how the orbits are interrelated. We have the
following differential fibrations with fibres $\CP^1$ for the
$\pi_0$ and the $\pi_2$ projections and $\CP^2$ for the $\pi_1$'s.
Again, as $U(1)=\tilde{U}(1)=SO(2)$, both $U(3)$ and $\tilde{U}(3)$
contain the above $U(2)\times U(1)$.

\begin{equation}\begin{array}{cccccc}

\begin{array}{c} \,\,\,\,\,\sO\\
\,\,\,\,\,\,\,\,\downarrow ^{\pi_0}\\ \,\,\,\,\,\,\,\sF^+\\

\begin{array}{ll}^{\pi_1}\swarrow & \,\,\, \searrow^{\pi_2}\\
 \sP^+ & \,\,\,\,\,\,\,\,\,\,\sG\end{array}
\end{array}

&\qquad &

\begin{array}{c} \sO\\
\,\,\,\,\downarrow ^{\pi'_0}\\ \sD\\

\begin{array}{ll} \,\,\,\,\,^{\pi^+_1}\swarrow & \,\,\, \searrow^{\pi^-_1}\\
 \,\sP^+ & \,\,\,\,\,\,\,\,\,\,\,\,\,\sP^-\end{array}
\end{array}

&\qquad &

\begin{array}{c} \sO\\
\,\,\,\,\downarrow ^{\pi''_0}\\ \sF^-\\

\begin{array}{ll} ^{\pi''_1}\swarrow & \,\,\, \searrow^{\pi''_2}\\
 \sG & \,\,\,\,\,\,\,\,\,\,\,\sP^-\end{array}
\end{array}
\end{array}\label{trianglefibr3}\end{equation}

\vskip10pt

\noindent Each fibre can be interpreted as a coadjoint orbit.

\begin{prop} (Bernatska-Holod \cite{BerHo}) Given a compact
  semisimple Lie group $G$, suppose that $G_\alpha$ (the isotropy
  group of $\alpha$) is not a maximal subgroup of $G$.  Then there
  exists a subgroup $H$ such that $G_\alpha\subset H\subset G$ and
  $\sO_\alpha$ fibres over $G/H$ with fibre
  $H/G_\alpha$.\label{bernat}\end{prop}

\noindent We express this in symbols by $\sO_\alpha\cong G/H\rtimes
H/G_\alpha$.\smallbreak

According to \cite{BerHo} a generic $SO(2N)$ orbit can be viewed as
follows:$$\sO^{SO(2N)}\cong \Gr_2(\R^{2N})\rtimes \sO^{SO(2N-2)},$$
which in our case becomes:
\begin{equation}\sO^{SO(6)}\cong\Gr_2(\R^{6})\rtimes \sO^{SO(4)}\cong
\Gr_2(\R^{6})\rtimes \Gr_2(\R^{4})\cong \Gr_2(\R^{6})\rtimes
(S^2\times S^2)\label{fibr.so6}\end{equation}

\noindent Furthermore observe that the only orbit of $SU(2)$ is
$\sO^{SU(2)}=\frac{SU(2)}{U(1)}\cong\CP^1$.\bigbreak

\noindent Then for  $SU(3)$, $\sO^{SU(3)}=\frac{SU(3)}{S(U(1)\times
U(1)\times U(1))}$ and $\sO^{SU(3)}_d=\frac{SU(3)}{S(U(2)\times
  U(1))}\cong\CP^2$.\medbreak

\noindent The generic $SU(3)$ orbit fibres over the degenerate one:
\begin{equation}\sO^{SU(3)}\cong\sO^{SU(3)}_d\rtimes\sO^{SU(2)}\cong\CP^2\rtimes\CP^1
\label{fibr.su3}\end{equation}\bigbreak

\noindent But $\sO^{SO(6)}$ fibres over $\sP^\pm$ with fibre
$\frac{U(3)}{U(1)\times(1)\times
U(1)}\cong\frac{SU(3)}{S(U(1)\times(1)\times U(1))}$ so:
\begin{equation}\sO^{SO(6)}\cong\sP\rtimes\sO^{SU(3)}
\cong\CP^3\rtimes\CP^2\rtimes\CP^1 \label{fibr.su4}\end{equation}

The above fibrations are \emph{symplectic} in the sense that their
fibre $\pi^{-1}(p)=F$ is a symplectic manifold for which the
transition mappings induce symplectomorphisms of $F$. This is
equivalent to a certain connection being flat \cite{GLSSFM}. More
generally, for any compact Lie group $G$, we have:

\begin{prop}(Guillemin-Lerman-Sternberg \cite{GLSSFM}) Let $x,\,\lambda$ be points in the same Weyl chamber.
  If the isotropy Lie algebras satisfy $\g_x\subset\g_\lambda$, then
  the map $\mathscr{O}_ x\to \mathscr{O}_\lambda$ given by $g\cdot
  x\mapsto g\cdot\lambda$ is a symplectic fibration with fibre a
  $G_\lambda$-coadjoint orbit.\label{cosimpfibr}
\end{prop}

In this context it is easy to prove that:

\begin{prop} A fixed point of a given orbit fibres over a fixed point in a lower one.
\end{prop}

\begin{figure}[!h]
\centering\includegraphics[width=0.3\textwidth]{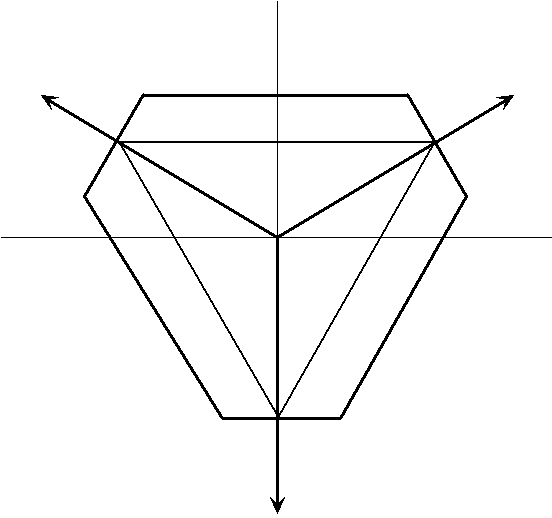}\\
\caption{A generic $SU(3)$ coadjoint orbit fibres symplectically
over $\CP^2$} \label{fibsu3}\end{figure}

\begin{figure}[!h]
\centering\includegraphics[width=1.0\textwidth]{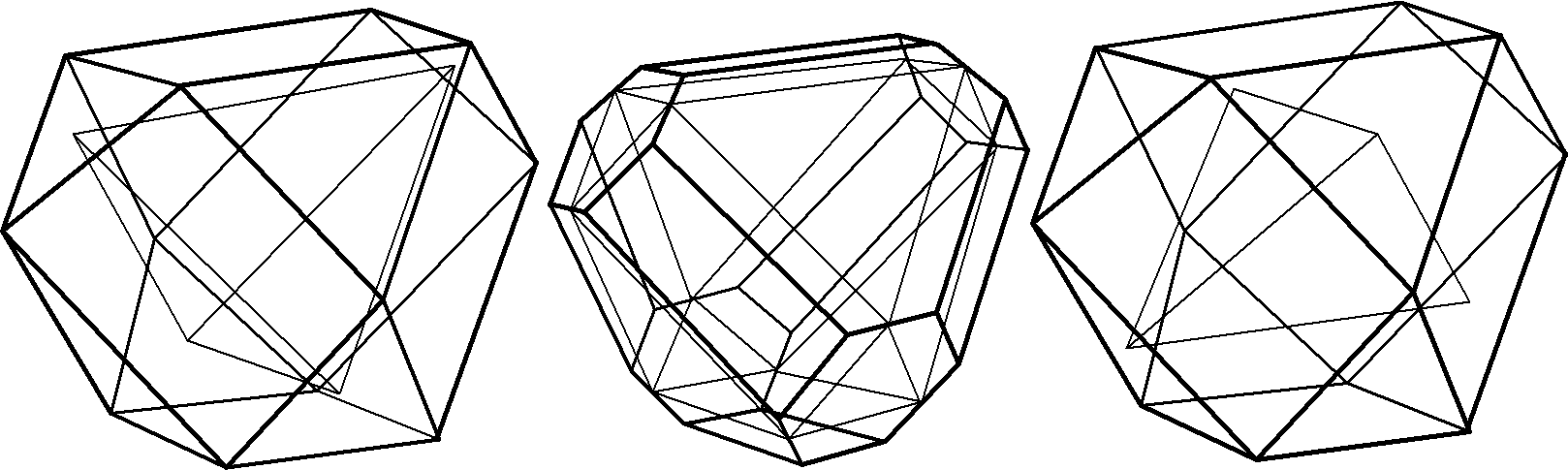}\\
\caption{$\sD$ fibres over $\sP^+$; \ $\sO^{SO(6)}$ fibres over over
$\sD$; \ $\sD$ fibres over $\sP^-$.} \label{fibsu3II}\end{figure}

\noindent This fact enables one to determine the image via $\mu_T$ of a fibre over fixed point in a lower orbit, i.e. the covex hull of fixed points contained in the fibre. Toric manifolds ($dimT=1/2\, dim M$) can be recognized by their moment (Delzant) polytopes. We will provide several examples of how symplectic fibrations over (symplectic submanifolds of) coadjoint orbits are efficiently illustrated by the moment map, even though the torus actions are typically low-dimensional and thus \emph{not} toric. For example the image of a generic $SU(3)$ orbit is a hexagon, and in Figure~\ref{fibsu3} we see how the symplectic fibration \rf{fibr.su3} is detected by the moment map. In fact, the $\CP^1$-fibres over the fixed points in $\CP^2$ (vertices of the triangle) are mapped to segments anchored at the fixed points of the generic orbit. Analogous considerations enable us to capture graphically the essence of the symplectic fibrations \rf{trianglefibr3} just comparing the moment polytopes as represented in  Figure~\ref{2fibrazioni} and Figure~\ref{fibsu3II}. A more detailed explanation of this observation emerges in the context of Section~3, other examples are given in Section~4.
 
\begin{rem}\label{diffweil}
  Observe that $U(2)\times U(1)$ is a subgroup of $SO(4)\times SO(2)$
  which suggests that $\sD$ fibres symplectically over $\sG$ with
  fibre $\CP^1$ (denote this projection by $\pi_2'$). In contrast with
  Proposition~\ref{cosimpfibr} we cannot find a representative of the
  lower orbit and an element of the fibre over it in the same Weyl
  chamber. For example compare the stabilizers of $(\alpha,\alpha,
  \beta)\in\sD$ and $(0,0,\beta)\in\sG$. We will return to this case
  in Section~4.
\end{rem}

\section{Riemannian geometry in six dimensions}

Let $G,H$ be Lie groups, $H$ a subgroup of $G$. It is a well known
fact that a reduction of a $G$-structure to an $H$-structure can be
realized by selecting a tensor $\xi$, stabilized by $H$ in a
suitable $G$-module. The parameter space of such reductions is the
$G$-orbit of $\xi$.\smallbreak

A section of the tensor bundle $T^\ast M\otimes TM$ can be regarded
as an endomorphism of each tangent space. Such an endomorphism can
be analysed in terms of its kernel and other eigenspaces, which (as
$p\in M$ varies) give rise to distributions. Integrability
properties of such distributions characterize the structure of $M$.

Let $(M,g)$ be a Riemannian manifold  of dimension $N$. Any smooth
$2$-form $\omega$ determines a skew-symmetric endomorphism $\fF$ of
each tangent space via
\begin{equation}
\omega(X,Y)=g(\fF X,Y). \label{endo}
\end{equation}

\noindent Suppose that $H$  is the largest non-trivial subgroup of
$SO(N)$ which preserves $\fF$ at every point in $M$. Then  $\fF$ determines a reduction of
the structure group $SO(N)$ to $H$. There are many relevant examples
of such a procedure. At one extreme (for N=2m) is the case when $\fF$ is an almost
complex structure on $T_pM$ ($H=U(m)$), the eigenspaces of $\fF$ are
maximal complex isotropic subspaces of $(T_pM)_c$. If a $2$-form is proportional 
to one associated to an almost complex structure, it \emph{determines} an almost
complex structure even though the corresponding $\fF$ is not an almost complex structure itself. 
In fact the eigenspaces of $\fF$ are again maximal complex isotropic subspaces of $(T_pM)_c$, and the eigenvalue is purely imaginary with respect to an underlying almost complex structure. 
At the other extreme if $\omega$ is a non-vanishing simple
2-form, then $(\ker\fF)^\perp$ has dimension 2,  and $\fF$ determines an
almost complex structure on this subspace  ($H=U(1)\times SO(N-2)$). In these
examples, $\fF$ satisfies $$g(X,Y)=g(\fF X,\fF Y),\qquad X,Y\in(\ker\fF)^\perp,$$ and thus
induces an orthogonal transformation on $(\ker\fF)^\perp$. For the
purposes of this article we say:\smallbreak

\begin{defi}\label{compy} Two distinguished Riemannian $G$-structures defined via
$2$-forms are called \emph{compatible} if the associated
skew-symmetric endomorphisms commute.
\end{defi}

The use of 2-forms to define a geometrical structure leads one
naturally to consider coadjoint orbits for $SO(N)$ which are complex
flag manifolds. An element of $\mathfrak{so}(N)$ can be regarded as
a skew-symmetric endomorphism. In even dimensions the eigenvalues of
such an endomorphism are pure-imaginary and paired and its spectral
structure is preserved by the orthogonal group. The flags in
question are determined by the set of eigenspaces of $\fF$ and the
stabilizer of the $2$-form depends only on the set of eigenspaces
and not on the precise eigenvalues. \smallbreak

In six dimensions
$\Lambda^2\mathbb{R}^6\cong\mathfrak{so}(6)\cong\mathbb{R}^{15}$.
The torus \rf{torus} acts in the standard block-diagonal way. The
images of the root spaces of $\theta_i$ are the subspaces $\langle
e^1,e^2\rangle$, $\langle e^3,e^4\rangle$, $\langle e^5,e^6\rangle$.
The Lie algebra $\mathfrak{t}$ is mapped to the $3$-dimensional
subspace of $\Lambda^2 T^*_pM$ spanned by the elements $\{ e^{12},
e^{34}, e^{56}\}$ (where $e^{ij}=e^i\wedge e^j$). Following this construction the image in
$\Lambda^2 T^*_pM$ of the fundamental Weyl chamber $B$ is generated
by  the elements
\begin{equation}
\label{mus} \mu_1=e^{12}+e^{34}+e^{56},\qquad\mu_2=e^{12},\qquad
\mu_3=e^{12}+e^{34}-e^{56}.
\end{equation}
\noindent Every coadjoint orbit has a unique representative element
in the closure $\bar{B}$ thus:

\begin{prop} Any $2$-form at a point of an oriented
  Riemannian 6-manifold is equivalent under the action of $SO(6)$ to
  a linear combination $\sum\limits_{i=1}^3a_i \mu_i$ with
  $a_i\ge0$.
\end{prop}

\medbreak

We can now introduce the Riemannian structure defined by a fixed
2-form. Such a structure is determined by a smooth section of the
fibre bundle $M\times_{SO(6)}\sO$, where $\sO$ is a coadjoint orbit.
The position of the representative of the orbit inside the image of
$\bar{B}$ determines the structure group of the reduction. Table 2, in which $a$, $a_i$ are positive,
relates positions inside $\bar{B}$ to $SO(6)$ orbits ordered by
dimension (recall Table 1).\bigbreak

\noindent\phantom.\hspace{10pt}
\begin{tabular}{|l|l|l|}\hline
                  Case & $2$-form                                                                     & SO(6) orbit\\ \hline\hline
                  1 & $a\mu_1=a(e^{12}+e^{34}+e^{56})$                                                & $ \sP^+$\\ \hline\hline
                  2 & $a\mu_3=a(e^{12}+e^{34}-e^{56})$                                                & $ \sP^-$\\ \hline\hline
                  3 & $a\mu_2=ae^{12}$                                                                & $ \sG$\\ \hline\hline
                  4 & $a_1\mu_1+a_2\mu_2=(a_1+a_2)e^{12}+a_1(e^{34}+e^{56})$                              & $\sF^+$\\ \hline\hline
                  5 & $a_1\mu_2+a_2\mu_3=(a_1+a_2)e^{56}+a_1(e^{34}-e^{56})$                                & $\sF^-$\\ \hline\hline
                  6 & $a_1\mu_1\+a_2\mu_3=(a_1\+a_2)(e^{12}\+e^{23})\+(a_1\-a_2)e^{56},\,\,\,\,\,a_1>a_2$ & $\sD^+$\\ \hline
                  7 & $a_1\mu_1\+a_1\mu_3=2a_1(e^{12}\+e^{34})$                                             & $\sD^0$\\ \hline
                  8 & $a_1\mu_1\+a_2\mu_3=(a_1\+a_2)(e^{12}\+e^{34})\+(a_1\-a_2)e^{56},\,\,\,\,\,a_1<a_2$ & $\sD^-$ \\ \hline
                  9 & $a_1\mu_1+a_2\mu_2+a_3\mu_3$                                                        & $\sO^{SO(6)}$\\ \hline
\end{tabular}\smallbreak

Table 2\bigbreak

\noindent We comment briefly case by case:\smallbreak

\noindent \textbf{Case 1}. The isotropy group is $U(3)$, so working
pointwise, without implying integrability, we shall refer to the
corresponding $G$-structure as an \emph{orthogonal almost complex
structures} (OCS) on $T_pM$ compatible with a fixed orientation. The
parameter space at each point of $M$ is $\sP^+$.\smallbreak

\noindent \textbf{Case 2}. $\sP^-$ parametrizes the OCS's inducing
the opposite orientation on $T_pM$.\smallbreak

\noindent \textbf{Case 3} features the Grassmannian $\Gr _2(\R^6)$
of oriented 2-planes in $\R^6$. A simple 2-form defines via
\rf{endo} a splitting $T_pM=\mathcal{V}\oplus\mathcal{H}$ with
$\mathcal V$ an oriented 2-plane, and $\mathcal H=\ker\fF$ a 4-plane
whose orientation is not specified. The orbit $\sG$ parametrizes a
set of \emph{orthogonal almost product structures} (OPS) studied by
Naveira in \cite{NavAPS} (alternatively defined by a
$(1,1)$-tensor field $P=v-h$ where $v$ and $h$ represent the
projections on $\mathcal{V}$ and on $\mathcal{H}$).\smallbreak

\noindent Cases 4-8 are related to the 10-dimensional
``intermediate'' complex flag manifold; in this case the
$\mathrm{Ad}$ action is characterized by two distinct pairs of
imaginary eigenvalues. The corresponding isotropy subgroup of
$SO(6)$ is isomorphic to $U(1)\times U(2)$.

\begin{defi} Let $M$ be an $N$-dimensional Riemannian manifold. A
  \emph{mixed structure} (MS) on $M$ is a reduction of the structure group
  to $U(p)\times U(q)$, where $2(p+q)=N$.  \label{mixed}
\end{defi}\noindent Such a structure is equivalent to the simultaneous
assignment of an OCS $J$ ($J^2=-I$) and an OPS $P$ ($P^2=I$) which
are compatible ($JP=PJ$). In our case $p=1$ and $q=2$ we set
$\mathcal{V}=\ker(P-I)$ so that $P$ is the identity on the
2-plane.\smallbreak

\noindent \textbf{Case 4}. The $2$-form belongs to the plane
generated by $\mu_1$ and $\mu_2$. The position of this point inside
the Weyl chamber $\bar B$ reflects the fact that $\omega$ is a
linear combination of a 2-form arising from an almost complex
structure $J$ and a simple 2-form arising from a positively-oriented
$J$-invariant 2-plane.

\begin{prop} A $\sF^+$ orbit parametrizes MS's
  determined by an OCS $J$ in $\sP^+$ and an OPS in $\sG$ whose
  2-plane is $J$-invariant and oriented consistently with $J$.
\end{prop}

\noindent \textbf{Case 5} is analogous; $\sF^-$ parametrizes MS's
with $J\in\sP^-$ and an OPS whose 2-plane is $J$-invariant and
oriented consistently with $J$.\smallbreak

\noindent \textbf{Cases 6 and 8}. This time, the position of the
2-form in $\bar B$ exhibits it as a weighted linear combination of
two compatible OCS's $J_+\in\sP^+$ and $J_-\in\sP^-$.

\begin{lemma} If two OCS's on $\R^6$ are compatible then they coincide
  up to sign on a real $4$-plane (and, therefore, on a complementary
  $2$-plane).\label{compplane}
\end{lemma}

\nit{Proof} Fix one of the structures, and use this to identify
$\R^6$ with $\C^3$. The second OCS lies in $U(3)$ with respect to
the first and its $3\times3$ matrix has eigenvalues $\pm i$, whilst
the first matrix is (say) $+i$ times the identity. In particular,
both matrices leave invariant a complex 2-dimensional subspace of
$\C^3$, giving rise to a real 4-plane.\qed\medbreak

\noindent As $J_+$ and $J_-$ belong to different orbits, they coincide
on an invariant $4$-plane, but induce opposite orientations on the
complementary $2$-plane. In the specific case displayed, the latter is
$\langle e^5,e^6\rangle$. Thus,

\begin{prop}\label{orbd} A $\sD^\pm$ orbit parametrizes MS's
  determined by an OCS $J$ in $\sP^\pm$ and an OPS in $\sG$ whose
  2-plane is oriented consistently with $-J$.
\end{prop}

\noindent \textbf{Case 7} is the special case in which the
contributions of the two OCS's have the same weight. The
2-dimensional subspace is determined by the kernel of $\omega$ and
thus the orientation is not specified by the tensor. This
corresponds exactly to Yano's definition of Riemannian $f$-structure
\cite{Yano, YKSM}, further developed by Blair \cite{BlGMG}:

\begin{defi} An $f$-structure on a differentiable manifold is a tensor $f$
(as the one in \rf{endo}) satisfying $f^3+f=0$, the existence of
which is equivalent to a reduction of the structure group to
$U(p)\times SO(q)$.\end{defi}

\noindent In conclusion, the $\sF$ and $\sD$ type orbits parametrize
at each point the $SO(6)$-inequivalent (but $O(6)$-equivalent) mixed
structures. The fact that $f$-structures in six dimensions provide a
special case of MS's is due to the isomorphism $SO(2)\cong U(1)$
(strictly speaking also the OPS's parametrized by $\sG$ are examples
of $f$-structures).\smallbreak

\noindent \textbf{Case 9} parametrizes the set of possible
$T^3$-reductions of the Riemannian structure, consisting of a choice
of three orthogonal complementary $2$-dimensional spaces in each
tangent space.\smallbreak

\begin{rem}
  Our construction refines the description of a geometrical structures,
  in the sense that we put more emphasis on the defining tensor rather
  than merely the isotropy subgroup. However the stabilizer of the $2$-form
  depends only on the spectral structure (set of eigenspaces) of the
  corresponding skew-symmetric endomorphism and does not depend  on
  the specific eigenvalues. Thus all the points staying
  in analogous positions in $\bar{B}$ represent
  \emph{the same} $G$-structure. The precise values of $a_i$ are not relevant.\end{rem}

  In the case of OCS's, OPS's and f-structures we can
  identify the $G$-structure with a real projective class in $\bar{B}$.
  All the points belonging to the same wall of $\bar{B}$ represent the same
  MS.

\section{Moment polytopes}

Previously we analysed a coadjoint orbit $\sO$ as a symplectic
manifold. The Riemannian $G$-structures under consideration are now
realized as smooth sections of fibre bundles with fibre $\sO$. The
mapping
$$\frac{SO(6)}{G}\to\Lambda^2T^*M$$
which associates a $2$-form to a specific $G$-reduction can be
interpreted (at each point of $M$) as the moment map
$$\frac{SO(6)}G\to\mathfrak{so}(6)^*$$
associated to the KKS symplectic structure. Combining this mapping
we obtain the orthogonal projection
$\mathfrak{so}^*(6)\to\mathfrak{t}^*$ which gives us the moment
mapping $$\mu_T:\frac{SO(6)}G
\longrightarrow\mathfrak{t}^*\cong\R^3$$ for the action of $T$
itself. To sum up,

\begin{teor} The Hamiltonian action of the maximum torus $T$ of
  $SO(6)$ on $\sO$ associates a characteristic ``moment
  polytope'' to each of the Riemannian structures defined by a
  2-form.\label{main}
\end{teor}

In Section 1 we proved, using standard Lie group theory, what the
precise shape of the moment polytope associated to each structure is
(compare Table~1 and Table~2). The aim of this section is to
describe how those polytopes can be obtained and interpreted in
terms of $2$-forms and compatibility of $G$-structures (in the sense
of Definition~\ref{compy}). For this purpose the subsets of $\sO$
consisting of points fixed under the action of a suitable subgroup
of the maximal torus $T=T^3$ provide in each orbit a ``skeleton" of
relevant structures. Those are the subsets on which $\mu_T$ is
singular.\smallbreak

Let $G$ be a compact Lie group. Having chosen a maximal torus, we
take a set of fundamental weights in $\g^*$. Denote by $\la_1,\ldots\la_N$ the set which includes the fundamental weights and all their conjugates, by $F_1,\ldots F_N$ the stabilizer group of each weight and by $W_i$ the Weyl group of the structure $(F_i,T)$.
Observe that $W_i$ is generated by the reflections induced by the
roots orthogonal to $\la_i$. As $F_i$ leaves invariant $\la_i$, we
have
\[\mathrm{ad}_X(\la_i )=[X,\la_i]=0,\qquad X \in \mathfrak{f}_i,\]
where $\mathfrak{f}_i$ is the Lie algebra of $F_i$.  The previous
equation can be read ``backwards'' to give
$\mathrm{ad}_{\la_i}(X)=0$, so the elements of subalgebra
$\mathfrak{f}_i\subset\g$ are preserved by the action of the circle
subgroup
$$C_i=\{\exp{2\pi\imath(t\la_i)}:t\in\R\}.$$\begin{teor}\label{important}
(Guillemin-Lerman-Sternberg \cite{GLSSFM})
  The critical (or singular) sets of the torus moment map $\mu_T:
  \mathscr{O_\lambda}\to\mathfrak{t}$ are the symplectic
  manifolds$$F_i\cdot w\lambda,\qquad w\in W,\quad i=1,\ldots,N.$$The critical
  values of $\mu_T$ are the corresponding convex polytopes $\mathrm{conv}(W_i\cdot w\lambda )$.
\end{teor}

An immediate application of this theorem to the case of $SO(6)$
allows one to visualize the critical values of $\mu_T$.

\begin{prop}Given a vertex $\alpha$ of the moment polytope $\Delta$ of an
$SO(6)$-coadjoint orbit, the image by $\mu_T$ of the symplectic
manifold $F_i\alpha$ consists of the intersection of $\Delta$ with
the plane orthogonal to $\la_i$, which passes through $\alpha$.
\label{weilweil}\end{prop}

\begin{figure}[!h]
\centering\includegraphics[width=1.0\textwidth]{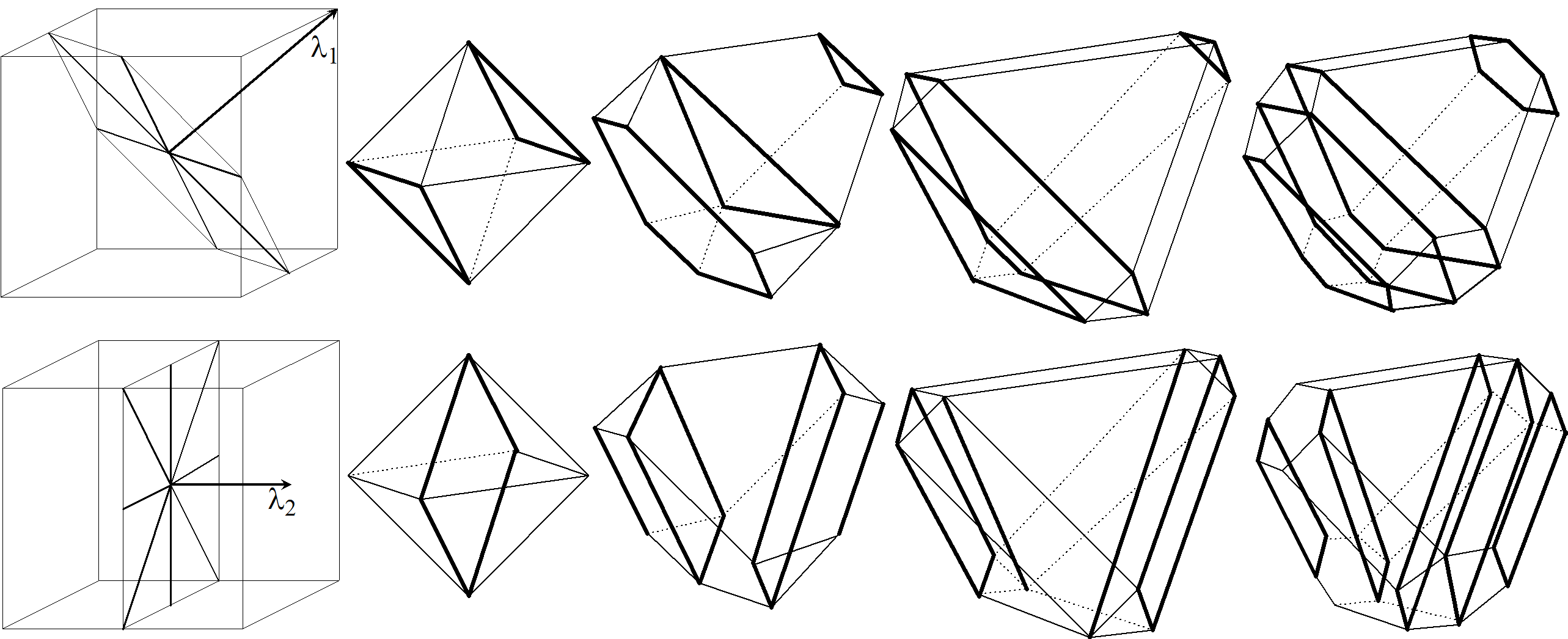}\\
\caption{Roots orthogonal to $\la_1$ and $\la_2$. Projections of
$C_1$ and $C_2$-invariant sets.}
 \label{refldef}\end{figure}

Figure~\ref{refldef} shows the directions of some roots orthogonal
to $\la_1$ and $\la_2$ and the corresponding fixed-point sets. The
roots generating $W_1$ and $W_2$ can be viewed as inward pointing
normal vectors of the polytopes $\mathrm{conv}(W_i\cdot
w\lambda)=\mu_T(F_i\cdot w\lambda)$.

\begin{rem}The sets of points fixed by $C_i$ are
  $F_i$-coadjoint orbits. The case of $SO(6)$ appears to be highly constrained
  in the sense that each critical set can be recognized by looking at its
  moment polytope. The fixed-point set of each $C_i$ carries an effective
  Hamiltonian action of a complementary subtorus in $T$ which is a maximal torus in $F_i$. The stabilizer of $\la_1$ is
  $SU(3)$, its coadjoint orbits have been described in Section~1.
  Since $\CP^2$ is a toric variety, the filled triangles on Fig.~\ref{refldef}
  identify fixed-point sets symplectomorphic to $\CP^2$. This implies that each
  filled hexagon represents necessarily a generic $SU(3)$-orbit. The stabilizer of
  $\la_2$ is $U(1)\times SO(4)$. The symplectic manifolds $\Gr_2(\R^4)\cong S^2\times S^2$ and
  $\CP^1$ are coadjoint orbits of $SO(4)$ and both are toric varieties, mapped
  respectively onto the filled rectangles and the bold line segments on Fig.~\ref{refldef}.
\label{subtor}
\end{rem}

The sets of points fixed simultaneously by the action of two 1-tori
are projected onto line segments determined as intersections of the
images of the single 1-tori fixed-point sets. Examples of such
projections are the edges of any polytope but also internal segments
and segments contained in faces. Certain segments appear as images
of sets invariant under $C_2$. This fact is quite easy to justify in
terms of sums of roots. \smallbreak

The Delzant Theorem \cite{delzant} states that the toric moment map defines a bijective correspondence between symplectic toric manifolds and Delzant polytopes. The following Proposition gives an operative criterion for establishing whether a symplectic submanifold of a $G$-coadjoint orbit is toric. The conditions on edges or vectors normal to faces, which define a Delzant polytope (simplicity, rationality and smoothness, see for example \cite{cannas}) are related to a lattice in $\mathfrak{t}^*$ determined by the root system of $G$.

\begin{prop}\label{refery} Let $G$ be a compact Lie group with maximal torus $T$, and let $\mathscr{R}\subset \mathfrak{t}^\ast$ be the set of roots. Let $M$ be a coadjoint orbit of $G$, and let $\mu_T:M\longrightarrow \mathfrak{t}^\ast$ be the moment map for the $T$-action. Let $N\subset M$ be a connected $T$-invariant symplectic submanifold and let $N^T\subset N$ denote the subset of points fixed by the torus action. Given a point $p\in N^T$ define:\medbreak

$\mathscr{R}_p=\{ \lambda\in \mathscr{R}|\mu_T(q)-\mu_T(p)=c\lambda$ for some $q\in N^T$ and $c>0\}$.\medbreak

\noindent If the vectors in $\mathscr{R}_p$ are linearly independent for some $p\in N^T$, then $N$ is a symplectic toric manifold.

\end{prop}

\nit{Proof} Since $N$ is connected, it is enough to prove that the weights for the $T$-action on $T_pN$ are linearly independent. Since $N\subset M$, the weights for the $T$ action on $T_pN$ are a subset of the weights for the $T$-action on $T_pM$, which themselves are a subset of $\mathscr{R}$. So we only need to consider weights in $\mathscr{R}$. Given $\lambda\in \mathscr{R}$, let $K\subset T$ be the kernel of the character associated to $\lambda$. Then the isotropy submanifold $N^K$ must contain at least one fxed point $q$ such that $\mu_T(q)-\mu_T(p)$ is a positive multiple of $\lambda$.\qed\medbreak

\noindent The statement and the proof of Proposition~\ref{refery} were suggested to the author by the reviewer of this article.\medbreak

A relation to the theory of G-structures issues from the following (which is obvious):
 
\begin{prop} If a $2$-form is fixed by the action of some subgroup $C$ of
the maximum torus $T\subset SO(N)$, then the corresponding
skew-symmetric endomorphism $\mathfrak{F}$ commutes with the action
of $C$ on $\R^N$.\end{prop}

\noindent More concretely an OCS acts as a simultaneous rotation by
$\pi/2$ on each of a triple of invariant $2$-planes and can be
interpreted as an element of a suitable one-torus. In particular,
the OCS associated to $\mu_1$ acts as
$\exp(\imath\frac{\pi}{2}\cdot\la_1)$ (see Figure~\ref{parameter2}).
Similarly, the endomorphism determined by $e^{12}$ corresponds to
$\exp(\imath\frac{\pi}{2}\cdot\la_2)$.\medbreak

We describe in this context the tetrahedra associated to both cases
of almost complex structures. This technique will be very useful
later on. The faces of each tetrahedron are projections of sets
fixed by the action of the circle generated by the ``opposite''
fundamental weight, so they represent OCS's commuting with the OCS
represented by the opposite vertex.

Lemma~\ref{compplane} implies that two commuting OCS's in the same
orbit coincide on a $2$-plane and differ by a sign on the
complementary $4$-plane. The sum of the corresponding $2$-forms
gives (twice) a simple $2$-form detecting the common invariant and
consistently-oriented $2$-plane. We uniform the notation to the one
of \cite{AGSHGI} setting $\mu_1=\omega_0$ (see \eqref{mus}):
$$\omega_0=e^{12}+e^{34}+e^{56},$$ \noindent and $J_0$ will denote the
corresponding OCS. The $2$-forms in the $SO(6)$ orbit of $\omega_0$,
which define OCS's commuting with $J_0$ are given
by:\begin{equation}\omega=-\omega_0+2v\wedge
J_0v,\label{formform}\end{equation} \noindent where
$v=\sum_{i=1}^6x_ie^i$ and $\sum_{i=1}^6x^2_i=1$. The prototypes are
the vertices of
$\Delta_{\sP^+}$:\begin{equation}\omega_1=+e^{12}-e^{34}-e^{56},\,\,\,\,\,\omega_2=-e^{12}+e^{34}-e^{56},\,\,\,\,\,
\omega_3=-e^{12}-e^{34}+e^{56}.\label{omegai}\end{equation}
\noindent   As $\mu_T$ is the projection to $\langle
e^{12},e^{34},e^{56}\rangle$ and
$$J_0v=x_1e^2-x_2e^1+x_3e^4-x_4e^3+x_5e^6-x_6e^5,$$
we have
\[\mu_T(\omega) =
(-1+2(x_1^2+x^2_2),-1+2(x^2_3+x^2_4),-1+2(x^2_5+x^2_6)) = (x,y,z).\]
Thus the projections of this form satisfy $x+y+z=-1$, and so lie in
the plane passing through $\omega_i$ with $i=1,2,3$. The planes
perpendicular to each $\omega_i$ and passing trough the
complementary vertices of the tetrahedron are obtained in the same
way. Applying  the same procedure to a generic $2$-plane generated
keeping the form inside $\sP^+$, the Cauchy-Schwarz inequality leads
to the condition $x+y+z\geq-1$ (see \cite{GMSSTV}). Varying the
vertex we get the entire set of inequalities that determine each of
the tetrahedra.\medbreak

Let us denote by $\Lambda^2_+ \R^4$ (respectively $\Lambda^2_-\R^4$)
the three-dimensional space of self-dual (anti self-dual) $2$-forms
on $\R^4$.

\begin{lemma}
Given an OCS $J$ on $\R^4$, the set of $J$-invariant and
consistently-oriented planes is isomorphic to $S^2$. \label{tec}
\end{lemma}

\noindent This is obvious because the set of complex lines in
$(\R^4,J)=\C^2$ is $\CP^1$. But to see the result in terms of
2-forms, recall that
$$\Lambda^2_+(\R^4)=\R\oplus\Lambda^{2,0}\oplus\Lambda^{0,2},\qquad
\Lambda^2_-(\R^4)=\Lambda_0^{1,1}.$$ A simple form can be written as
a sum of elements in $\Lambda^2_+(\R^4)$ and $\Lambda^2_-(\R^4)$ of
equal norm, and a $J$-invariant simple $2$-form is given by the
expression $$v\wedge Jv=\eta_++\eta_-,$$ \noindent where
$\eta_\pm\in\Lambda^2_\pm(\R^4)$ and $|\eta_+|=|\eta_-|$.\medbreak

Lemma~\ref{compplane} and Lemma~\ref{tec} imply the following facts:

\begin{corol}\label{important1} Given an OCS $J\in\sP^+$, the
subset of $\sP^\pm$ of OCS's compatible with $J$, is in one-to-one
correspondence with the set of $J$-invariant $2$-planes consistently
oriented with $\pm J$.
\end{corol}

\noindent Given an OCS on $\R^6$ compatible with $J$, there is an
$S^2$ of common invariant planes oriented consistently with $-J$
inside the common invariant $4$-plane. It is easy to prove:

\begin{corol}\label{important2} Given an OCS $J\in \sP^+$ and an
$J$-invariant $2$-plane $\alpha$ oriented consistently with $-J$,
there is an $S^2\in \sP^+$ of OCS's compatible with $J$, for which
$\alpha$ is invariant and consistently oriented.
\end{corol}

The edges of the tetrahedra (intersections of two faces) are
projections of OCS's commuting simultaneously with the two vertices
opposite to the faces in question. Lemma~\ref{compplane} and
Lemma~\ref{tec} imply that these sets are $2$-spheres. For example
the edges which do not contain $\omega_0$ are projections of forms
with fixed norm in:\begin{equation}\label{edges}
-e^{12}+\Lambda^2_-\langle e^3,e^4,e^5,e^6\rangle,\,\,\,\,\,\,
-e^{34}+\Lambda^2_-\langle e^1,e^2,e^5,e^6\rangle,\,\,\,\,\,\,
-e^{56}+\Lambda^2_-\langle e^1,e^2,e^3,e^4\rangle \end{equation}

\begin{figure} \centering\includegraphics[width=0.4\textwidth]{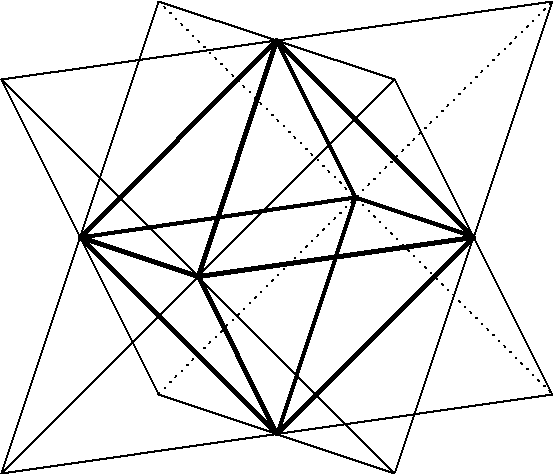}\\
\caption{The moment polytope $\Delta_\mathcal{G}$ interpreted as
intersection of $J$-invariant sets.} \label{planesI}\end{figure}

In this context we can also determine the faces of the octahedron
$\Delta_\sG$. Let $J$ be a vertex of $\Delta_{\sP^+}$ or
$\Delta_{\sP^-}$; we can ask ``what is the image in $\Delta_\sG$ of
the $J$-invariant planes''. The octahedron combines the inequalities
defining $\Delta_{\sP^+}$ and $\Delta_{\sP^-}$;$$|x|+|y|+|z|\le1,$$
\noindent thus $\Delta_{\sG}$ is obtained as an intersection of the
tetrahedra as in Figure~\ref{planesI}. This resumes the argument
applied in \cite{GMSSTV} to determine $\mu_T(\sG)$. In the next
section we will provide analogous interpretation of the moment
polytopes related to the remaining $SO(6)$ coadjoint orbits.

\section{A Klein correspondence}

In Section 2, we introduced the symplectic fibrations of $SO(6)$
coadjoint orbits. Consider the lower part of the first and the third
diagram in \rf{trianglefibr3}. The projections $\pi_1$ and $\pi_2$
(resp. $\pi_1''$, $\pi_2''$) can be understood in terms of the
classical Klein correspondence in which
$\sG=\Gr_2(\R^6)\cong\Gr_2(\C^4)$ is identified with a
non-degenerate quadric in $\mathbb{P}(\Lambda^2\C^4)$:\smallbreak

1. $\sG$ parametrizes the projective lines $\CP^1$ in
$\CP^3$.\smallbreak

2. A point $x\in\CP^3$ determines an $\alpha$-plane in $\sG$,
consisting of all the lines passing through that point.\smallbreak

3. A point $y\in(\CP^3)^*$ determines a $\beta$-plane in $\sG$,
consisting of all the lines lying in the plane $y$.\smallbreak

\noindent In the light of our realization of coadjoint orbits as
parameter spaces of Riemannian structures, this correspondence
assumes a completely new interpretation. Namely,\smallbreak

$1'.$ Given a decomposition $T_pM=\mathcal{V}\oplus\mathcal{H}$
arising from an OPS $P$, there is a $\CP^1$ worth of compatible
OCS's parametrized by $\omega\in S^2\subset
\Lambda^2_+\mathcal{H}^*$. This is our projective line in
$\sP^+$.\smallbreak

$2'.$ Given an OCS $J$ we have the $J$-invariant $2$-planes
generated by $\{v,Jv\}$ and each one determines an OPS.\smallbreak

$3'.$ Likewise, given an OCS $J$ we have the $J$-invariant
oppositely-oriented $2$-planes generated by $\{v,-Jv\}$.\smallbreak

\noindent To understand $1'$, recall that the 2-sphere of unit
self-dual forms parametrizes OCS's on $\mathcal H=\R^4$ compatible
with both metric and orientation. When combined with a standard
almost complex structure on $\mathcal V$, we obtain a
positively-oriented OCS on $\R^6$.\smallbreak

The results of Section~3 can be exploited to establish a mapping
between moment polytopes induced by the Klein correspondence. A MS
in $\sF^+$ determines an OCS $J$, namely its projection via $\pi_1$.
This $J$ identifies the tangent space $T_pM$ with $\C^3$, and
$J$-invariant splittings of $\R^6$ are parametrized by complex
1-dimensional (or complementary 2-dimensional) subspaces in $\C^3$,
i.e.\ by the projective space $\CP^2$. Thus, the set of MS's fibres
over the set of compatible OCS's with fibre $\CP^2$.

\begin{figure}[!h]

 \centering\includegraphics[width=0.8\textwidth]{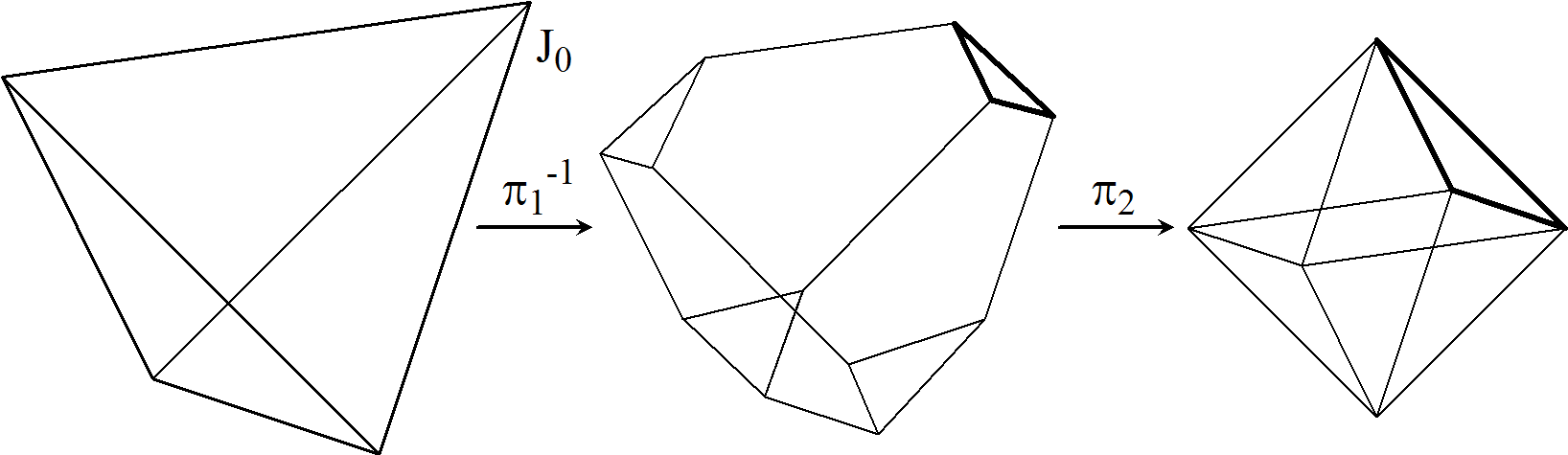}\\
\caption{A polytope map induced by the Klein correspondence.}
 \label{klein1}\end{figure}

The inverse image by $\pi_1$ of an OCS $J\in\sP^+$ inside $\sF^+$ is
determined by answering the question: ``which are the OPS's whose
2-plane is both $J$-invariant and oriented consistently with $J$?''
Working in terms of 2-forms, we take the non-degenerate 2-form
$\omega$ associated to $J$ and add to it a simple $2$-form $v\wedge
Jv$ representing the $J$-invariant plane in question. For instance,
$\mu_T(\pi^{-1}_1(J_0))$ can be easily determined by the technique
introduced in Section~3, we see that $\mu_T(\omega_0+\alpha (v\wedge
J_0v))$ belongs to the
plane\begin{equation}x+y+z=3+\alpha.\end{equation}\noindent In view
of the results of Section~2 and the symplectic nature of the
fibrations (recall Proposition \ref{main}), we conclude that the set
of MS's compatible with $J_0$ gets mapped by $\mu_T$ onto the
triangular face of the truncated tetrahedron, generated by the
vertices \begin{equation}\label{vertices}e^{12}+ e^{34}+ e^{56}+
\alpha e^{12},\quad e^{12}+ e^{34}+ e^{56}+ \alpha e^{34},\quad
e^{12}+ e^{34}+ e^{56}+ \alpha e^{56}.\end{equation} \noindent This
is a subset of $\sF^+$ of points fixed by
$\exp(\imath\frac{\pi}{2}\cdot\la_1)$. Lemma~\ref{tec} gives the
possibility to interpret the edges of this triangular face. For
instance the edge joining the first two vertices in \eqref{vertices}
is the projection of a set invariant under $C_1$ and $C_2$. It
represents the $S^2$ of $J_0$-invariant and consistently oriented
planes in the span $\langle e^1,e^2, e^3, e^4\rangle$.
Figure~\ref{klein1} represents the map between polytopes induced by
the maps $\pi_1$ and $\pi_2$.

We can now test the mapping induced between polytopes on subsets of
the lower orbits in the symplectic fibrations \rf{trianglefibr3}.
Let $J_i$ be the OCS's that correspond to the 2-forms $\omega_i$ in
\rf{omegai} i.e. the vertices in the tetrahedron $\mu_T(\sP^+)$.
Adopting notation from \cite{AGSAHG}, we call $\mathcal{E}_{12}$ the
edge joining the vertices $\omega_1$ and $\omega_2$ as in
Figure~\ref{invsqr}.

\begin{prop} Denote by $L:=\mu_T^{-1}(\mathcal{E}_{12})$. The image
  $\mu_T(\pi_1^{-1}(L))\subset\Delta_{\mathscr{F^+}}$ is the polytope shown on the
  left in Figure~\ref{invsqr}. The set $\pi_1^{-1}(L)$ is a symplectic toric manifold.
\label{edge}\end{prop}

\nit{Proof} The coadjoint orbits of any compact classical Lie group
admit a standard invariant complex structure $J$ defined by the
action of a maximum torus on the root spaces of the isotropy
representation (see \cite{ABEM}). The KKS form is the K\"ahler form
of a canonical K\"ahler structure compatible with $J$. The
projections $\pi$ intertwine the complex structure. Therefore the
preimage of a complex submanifold (observe that $L\cong\CP^1$) is a
complex submanifold, and so it is a symplectic submanifold. For this
reason, its moment image is the convex hull of its fixed points.
Hence, one just needs to determine the fixed points which map to the
desired set. But this is trivial as fixed points in the higher orbit
belong to fibres over fixed points in the lower one. Now $\mu_T(\pi_1^{-1}(L))$ is $3$-dimensional and $\pi_1^{-1}(L)$ is a real $6$-dimensional symplectic manifold (a $\CP^2$ bundle over $\CP^1$), so $\pi_1^{-1}(L)$ is $T$-invariant. Looking at $\mu_T(\pi_1^{-1}(L))$, it is easy to check that  the hypotheses of Proposition~\ref{refery} are satisfied so $\pi_1^{-1}(L)$ is toric.\qed\smallbreak

\begin{rem} Proposition~\ref{edge} confirms that the mapping between polytopes arises
  from symplectic phenomena. Also in this case the moment map captures the essence of
  the symplectic fibration. The inverse image of each point of $\sP^+$
  in $\sF^+$ is isomorphic to a $\CP^2$, and the edge of the
  tetrahedron is the projection of a $\CP^1$. The inverse image of the
  entire set is thus a symplectic $\CP^2$-bundle over $\CP^1$ and the subset $\pi_1^{-1}(L)$ projects to ``a
  triangle times a line''!
\label{edgeprism1}\end{rem}\medbreak

Analogously, the inverse image by $\pi_2^{-1}$ of an element in
$\sG$ is parametrized by suitably-oriented OCS's on $\mathcal H$.
Given an orthonormal basis $\{f_1,f_2\}$ of $\mathcal V$, we extend
the $2$-form $f^{12}$ by adding a unit element of
$\Lambda^2_+\mathcal{H}$ or inside $\Lambda^2_-\mathcal{H}$ so as to
obtain an OCS on $T_pM$. For example, $e^{56}$ can be completed to
$\omega_0$ and $\omega_3$ in $\sP^+$ (recall \rf{omegai}) and a
whole 2-sphere of similar non-degenerate 2-forms. The same simple
2-form can be also completed to $-\omega_1$ and $-\omega_2$ in
$\sP^-$ and an $S^2$ worth of OCS's with the opposite orientation.\smallbreak

The proof of Proposition~\ref{edge} holds in the following two cases.  

\begin{prop} \label{square} Denote by $K$ the set of $2$-planes in
  the subspace $\langle e^1,e^2,e^3,e^4\rangle$. The image
  $\mu_T(\pi_2^{-1}(K))$ is the rectangular prismoid represented in
  bold on the right of Figure~\ref{invsqr}. The set $\pi_2^{-1}(K)$ is a symplectic toric manifold.

\begin{figure}[!h]

\centering\includegraphics[width=1.0\textwidth]{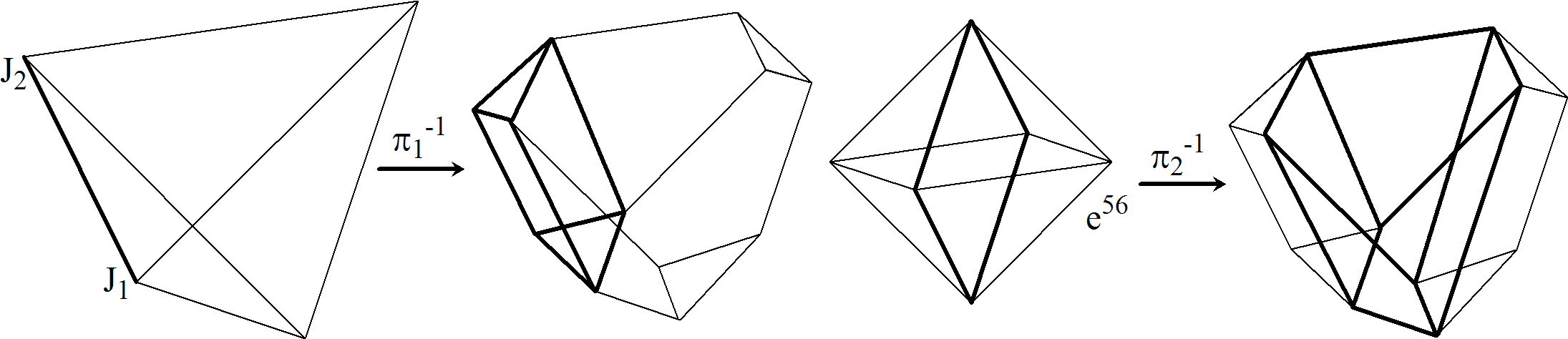}\\
\caption{Projections of the inverse images $\mu_T(\pi^{-1}_1(L))$
and $\mu_T(\pi^{-1}_2(K))$}
 \label{invsqr}\end{figure}
\end{prop}

\begin{rem} The image of $K\cong \Gr_2(\R^4)$ by the moment map is the square intersection
  of $\Delta_{\sG}$ with the plane $\langle e^{12},e^{34}\rangle$. In the
  present context it should be interpreted as a projection of
  $S^2\times S^2$ having $2$-spheres respectively in
  $\Lambda^2_+\langle e^1,e^2,e^3,e^4\rangle$ and $\Lambda^2_-\langle
  e^1,e^2,e^3,e^4\rangle$. We expect the subset
  $\pi_2^{-1}(K)$ of $\sF^+$ to be a symplectic $\CP^1$ bundle over $S^2\times S^2$.
  The intersection of $\mu_T(\pi_2^{-1}(K))$ with any plane orthogonal to $\langle e^{56}\rangle$
  is a rectangle, so this moment polytope is ``a rectangle times a line''.  \label{edgeprism2}\end{rem}\bigbreak

\begin{figure}[!h]
\centering\includegraphics[width=1.0\textwidth]{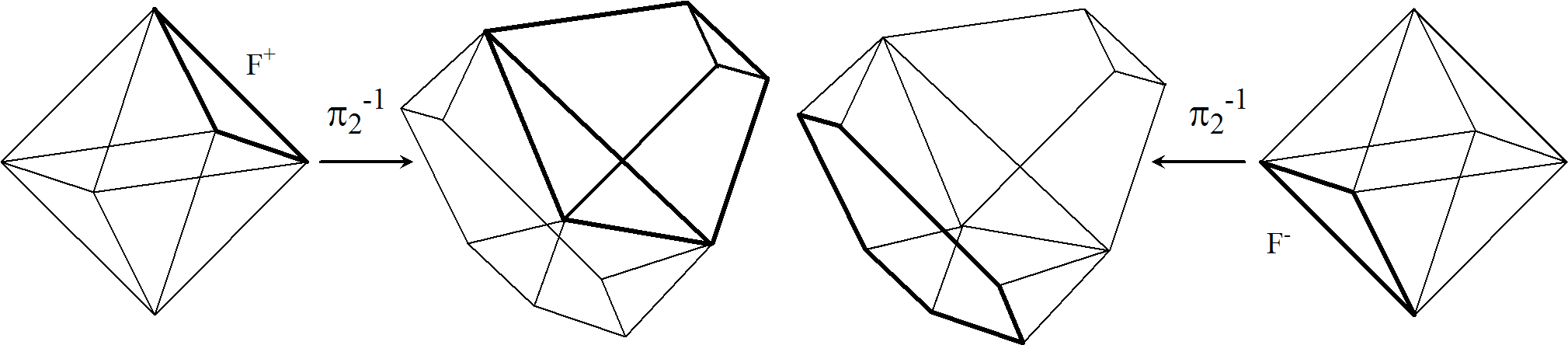}\\
\caption{ Projections $\mu_T(F^\pm)$ in $\Delta_{\mathscr{G}}$ and
$\mu_T(\pi_2^{-1}(F^\pm))$ in $\Delta_{\mathscr{F^+}}$.}
\label{newex}\end{figure}

\begin{prop} \label{newprop}Let $F^+$ and $F^-$ denote the $C_1$-invariant
  subsets of $\mathscr{G}$ (both symplectomorphic to $\CP^2$) projected by $\mu_T$
  on two disjoint faces of $\Delta_\mathscr{G}$. The images $\mu_T(\pi_2^{-1}(F^\pm))$
  are respectively a hexagon and a triangular prismoid as shown in Figure~\ref{newex}. 
  The set $\pi_2^{-1}(F^+)$ is a symplectic toric manifold.
\end{prop}

The sets $F^\pm$ represent $J_0$-invariant planes. Both
$\pi_2^{-1}(F^\pm)$ are symplectic fibrations over $\CP^2$ with
fibre $\CP^1$. The set $\pi_2^{-1}(F^\pm)$ contains the MS's defined
by $J_0$ itself (a $\CP^2$ projecting on the triangular face of
$\Delta_{\sF^+}$ on Fig.~\ref{newex}). The set of MS's obtained by a
plane in $F^+$ and the unique (by Corollary~\ref{important1})
corresponding OCS compatible with $J_0$ is a $\CP^2$ of MS's
compatible with $J_0$, which projects on the internal triangle.
Corollary~\ref{important2} implies the fibre over each point of
$F^-$ represents MS's compatible with $J_0$. This set projects on a
hexagonal face of $\Delta_{\sF^+}$. Alternatively:

\begin{rem}\label{hexa} The hexagonal face of
$\Delta_{\sF^+}$ represent MS's determined by an OCS compatible with
$J_0$ and a $J_0$-invariant plane oriented consistently with $-J_0$.
MS's compatible with $J_0$ are parametrized by $SU(3)$ orbits, and
the generic orbit is a symplectic $\CP^1$ bundle over
$\CP^2$.\end{rem}\smallbreak

The inverse image in $\sD$ of a point $\omega\in\sP^\pm$ is given by
all the forms obtained by adding to $\omega$ a ``small''
contribution negatively-oriented $J$-invariant $2$-plane. The
inverse image of a $J$ in $\sD$ can be recovered replacing the
positively-oriented $2$-plane of the case $\sF^+$ with a
negatively-oriented one. The map $(\pi^+_1)^{-1}\pi^-_1$ sends any
OCS in $\sP^+$ in the corresponding $\CP^2$ of compatible OCS's in
$\sP^+$ and vice versa. This concludes our ``$2$-form
interpretation'' of the middle symplectic fibration in
\rf{trianglefibr3}.\smallbreak

\begin{figure}[!h]
 \centering\includegraphics[width=1.0\textwidth]{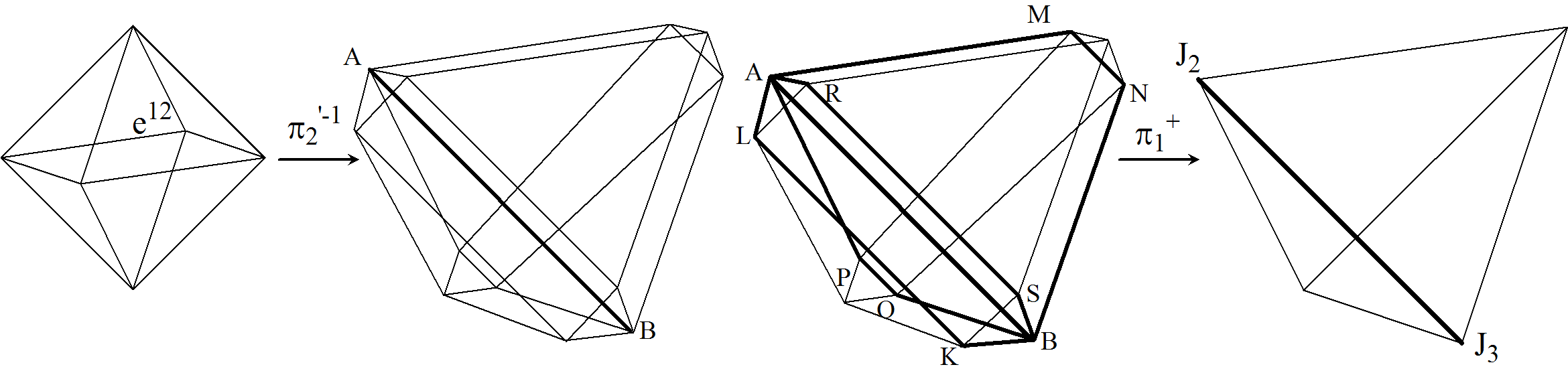}\\
\caption{A polytope map induced by the Klein correspondence.}
 \label{klein3}\end{figure}

In Remark \ref{diffweil} we observed that $\sD$ should fibre over
$\sG$. This fact becomes obvious in view of Proposition \ref{orbd}.
We can determine the inverse image of the point $e^{12}\in\sG$. This
case requires a 2-plane that is invariant simultaneously by $J_0$
and $-J_1$. The subsets in $\sD$ invariant by the rotations
generated by the roots corresponding to $J_0\in\sP^+$ and
$-J_2\in\sP^-$ project respectively onto the hexagons $ARSBOP$ and
$AMNBKL$ in Figure~\ref{klein3}. The projection of the inverse image
of $e^{12}$ is the internal segment $AB$, remarkably it is obtained
as the intersection of the above invariant sets. Now comparing
Figure~\ref{klein3} and Figure~\ref{parameter2} we see that there is
no representative of the fibre in the Weyl chamber where $e^{12}$
stays.

The image of $e^{12}$ in $\sP^+$ consists of all the OCS's $J$ for
which this form represents the contribution of a negatively-oriented
$J$-invariant $2$-plane (for example, full-rank forms involving
$-e^{12}$). It is easy to see that this set projects onto the edge
$\mathcal{E}_{23}$. Keeping in mind Lemma~\ref{compplane} it should
be not a surprise that the image of $F^+$ in $\sP^+$ projects onto
the face in $\Delta_{\sP^+}$ opposite to $J_0$.

As another example, consider the form $e^{15}-e^{26}-ae^{34}$. Its
image by $\mu_T$ is the midpoint of the segment $AB$.  It can be
interpreted as the sum of $e^{15}-e^{26}-e^{34}$ in $\sP^+$ and
$(1-a)e^{34}$. The form $e^{15}-e^{26}-e^{34}$ maps to the midpoint
of the edge $\mathcal{E}_{13}$.\smallbreak

\begin{figure}[!h]
\centering\includegraphics[width=1.0\textwidth]{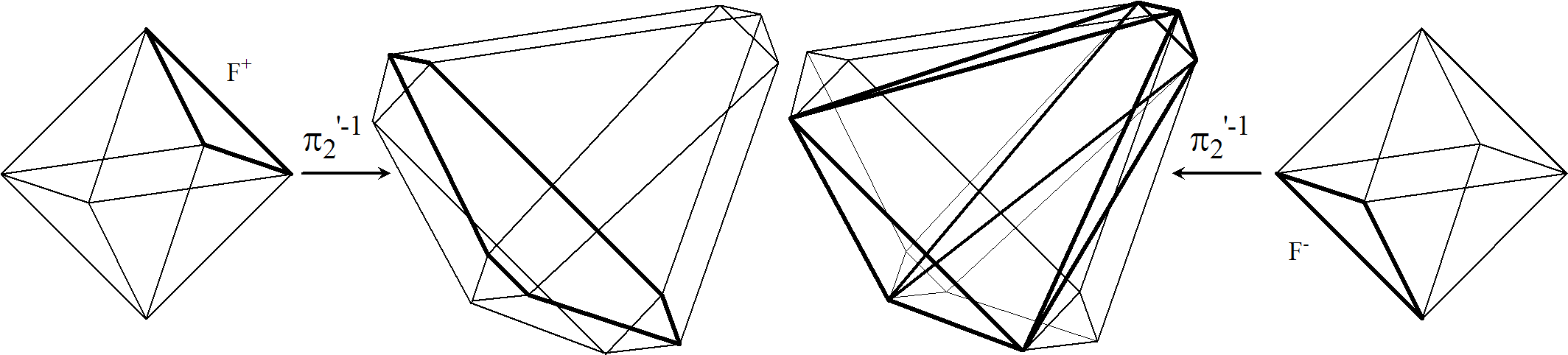}\\
\caption{Projections $\mu_T(F^\pm)$ in $\Delta_{\mathscr{G}}$ and
$\mu_T(\pi_2'^{-1}(F^\pm))$ in $\Delta_{\mathscr{D}}$.}
\label{newexample}\end{figure}

We invite the reader to complete the analysis of the maps between
the moment polytopes $\Delta_{\sP^\pm}$, $\Delta_{\sG}$ and
$\Delta_{\sD^\pm}$ keeping in mind Proposition~\ref{orbd} and
Remark~\ref{hexa}. In particular analogous considerations enable one
reader to interpret Figure~\ref{newexample}.

The moment polytope corresponding to $\sD^0$
(Figure~\ref{parameter1}g) is a special case of $\Delta_{\sD}$. The
rectangular faces are squares. For the $f$-structures, there is no
way to define a map between moment polytopes corresponding to the
symplectic fibrations in the middle diagram of \rf{trianglefibr3}.
The lack of such a mapping is due to the degeneracy of the
characteristic $2$-form. In fact, in the case of the remaining $\sD$
orbits, there are two distinguished internal segments corresponding
to the projections of the same $2$-plane taken with two different
orientations. In this case, the $2$-form is not capable of
determining the compatible orientation on its kernel. Graphically,
this is expressed by the fact that the two segments (corresponding
to the different orientations of the $2$-dimensional kernel)
intersect in the origin (see Figure~\ref{parameter1}(g) ). For
example the form $e^{14}+e^{23}\in\sD^0$ determines a splitting $
\langle e^1,e^2,e^3,e^4\rangle\oplus \langle e^5,e^6\rangle$ and can be extended to $e^{14}+e^{23}+e^{56}$ to yield a
compatible OCS in $\sP^+$. The image of $e^{14}+e^{23}$ by $\mu_T$
is the origin of $\R^3$, whereas the image of $\omega$ is the
midpoint $(0,0,1)$ of an edge of the tetrahedron.

\section{An application}

The intrinsic torsion of a geometrical structure is the first order
obstruction to its integrability. For this reason, a standard way of
classifying Riemannian $G$-structures is based on criteria whereby its
intrinsic torsion tensor $\tau$ reduces to a specific subset of
$G$-irreducible components of the corresponding space of intrinsic
torsion $\mathscr{W}$. For Riemannian $G$-structures, $\tau$ is
determined by the Levi-Civita derivative of the defining tensor and
$\mathscr{W}$ is isomorphic to $T^*M\otimes\mathfrak{g}^\bot$, where
$\mathfrak{g}^\bot$ is the orthogonal complement of the Lie algebra of
$G$ in $\mathfrak{so}(N)$. The constraint can be stated requiring that
some components of $\tau$ vanish, so we adopt the terminology
\emph{null-torsion classes}.\smallbreak

The prototype case gave rise to the sixteen classes of almost
Hermitian manifolds \`a la Gray--Hervella \cite{GrHSCA}. The
$U(n)$-irreducible components $\mathscr{W}_i\subset \mathscr{W}$ have
also been described in \cite{AHSSDF} by complexifying the exterior
algebra; the space $\Lambda^{p,q}\oplus\Lambda^{q,p}$ is the
complexification of a real vector space that is denoted
$[\![\Lambda^{p,q}]\!]$. Denote by $R(\lambda)$ the irreducible
complex $U(n)$-representations with dominant weight $\lambda$, and by
$\Lambda_0^{p,q}$ the Hermitian complement of the image of
$\Lambda^{p-1,q-1}$ under wedging with $\omega$.
We then have the following isomorphisms:
\begin{equation}\label{GH}
\begin{array}{c}
\mathscr{W}\cong T^*M\otimes\mathfrak{u}^\bot\cong\Lambda^{1,0}\otimes
\rr{\Lambda^{2,0}}\cong\mathscr{W}_1\oplus\mathscr{W}_2\oplus
\mathscr{W}_1\oplus\mathscr{W}_4,\\[5pt]
\mathscr{W}_1\cong \rr{\Lambda^{3,0}},\quad
\mathscr{W}_2\cong\rr{R(2,1,0,...,0)},\quad
\mathscr{W}_3\cong\rr{\Lambda_0^{2,1}},\quad
\mathscr{W}_4\cong\rr{\Lambda^{1,0}}.
\end{array}\end{equation}
The irreducible components immediately above give rise to well known classes of
almost Hermitian structures; $\tau$ lying in $\mathcal{W}_1$ means
that the structure is \emph{nearly-K\"ahler}, $\mathcal{W}_2$ --
\emph{almost-K\"ahler}, $\mathcal{W}_3$ -- \emph{cosymplectic
  Hermitian}, $\mathcal{W}_4$ -- \emph{locally conformal K\"ahler}
etc. In particular, the component of $\tau$ in $\mathcal{W}_1\oplus
\mathcal{W}_2$ can be identified with the Nijenhuis tensor, and so
\emph{Hermitian} structures belong to the class $\mathcal{W}_3\oplus
\mathcal{W}_4$.\medbreak

Analogous classification of the OPS's, developed by Naveira in
\cite{NavAPS} exploits the decomposition of the intrinsic torsion
space
\[\mathscr{V}\cong T^*M\otimes
\big(\mathfrak{so}(\mathcal{V})\oplus\mathfrak{so}(\mathcal{H})\big)
\cong(\mathcal{H}\oplus\mathcal{V})\otimes
\big(\mathcal{H}\otimes\mathcal{V}\big)\]
into the irreducible components
\begin{equation}
\label{VVVVVV}
\begin{array}{lll} \mathscr{V}_1=\Lambda^2
  \mathcal{V}\otimes\mathcal{H},\quad &\mathscr{V}_2=S_0^2
  \mathcal{V}\otimes \mathcal{H},\quad &
  \mathscr{V}_3=1_\mathcal{V}\otimes
  \mathcal{H},\\
  \mathscr{V}_4=\Lambda^2 \mathcal{H}\otimes\mathcal{V},
  &\mathscr{V}_5=S_0^2 \mathcal{H}\otimes\mathcal{V},
  &\mathscr{V}_6=1_\mathcal{H}\otimes\mathcal{V},
\end{array}
\end{equation}
using the notation of Section 2.

We denote by $\mathscr{W}_{ij\ldots}$ and $\mathscr{V}_{ij\ldots}$
the spaces $\mathscr{W}_i\oplus\mathscr{W}_j\oplus\ldots$ and
$\mathscr{V}_i\oplus\mathscr{V}_j\oplus\ldots$, or the corresponding
null-torsion classes.

In \cite{GMMSRS,GMSSTV} the same approach has been applied to MS's.
In six dimensions we exploit the geometrical interpretation of a
$U(1)\times U(2)$-structure in terms of the underlying OPS and OCS.
This allows the null-torsion classes of these MS's to be described by
means of the $U(1)\times U(2)$ module
\begin{equation}\label{uu}
\mathscr{M} \cong
T^*\!M\otimes(\mathfrak{u}(1)\oplus\mathfrak{u}(2))^\bot\cong
(\mH\oplus\mV)\otimes
\Big((\mV\otimes\mH)\oplus\rr{\la^{2,0}}\Big),
\end{equation}
where $\mH=\rr{\la^{1,0}}$ in accordance with the notation introduced
for \eqref{GH}. We shall also write $\mV=\rr{\nu}$, so that $\nu$ is a
complex vector space of dimension 1 on which $U(1)$ acts, and
$\Lambda^{1,0}=\nu\oplus\la^{0,1}$. Then $\mathscr{M}$ can be shown to
be isomorphic to the direct sum of 16 irreducible real summands:
\begin{equation}\label{fine}
\begin{array}{l}
\mathscr{M}\ \cong\
3\mH\oplus
2\mV\oplus
2\rr{\nu\la^{1,1}_0}\oplus
\rr{\nu^2\la^{1,0}}\oplus
\rr{\nu^2\la^{0,1}}
\\\hskip30pt\oplus\>
2\rr{\nu\la^{2,0}}\oplus
2\rr{\nu\la^{0,2}}\oplus
\rr{\nu\sigma^{2,0}}\oplus
\rr{\nu\sigma^{0,2}}\oplus
\rr{R(2,1)}.
\end{array}
\end{equation}
(Tensor product signs are omitted, $3\mH$ means
$\mH\oplus\mH\oplus\mH$. Also $\nu^2=\otimes^2\nu$, and
$\sigma^{2,0}$ is the second symmetric power of $\la^{1,0}$.) The 10
non-isotypic summands have respective dimensions
4,2,6,2,2,2,2,6,6,2. This approach enables one to compare the
intrinsic torsion of interrelated structures. For example, the fact
that\begin{equation}\label{MWV} \mathscr{M}=\mathscr{W} +
\mathscr{V}
\end{equation} implies the next result, also proved in
\cite{GMSSTV}.

\begin{prop}\label{nablajp} The intrinsic torsion tensor
  $\tau_\mathscr{M}$ of a MS is completely determined by the intrinsic
  torsion tensors $\tau_\mathscr{W},\tau_\mathscr{V}$ of the
  underlying OCS and OPS. Conversely, $\tau_\mathscr{M}$ determines
  the pair $(\tau_\mathscr{W},\tau_\mathscr{V})$.
\end{prop}

\noindent Expression \eqref{MWV} is not a direct sum as
$\tau_\mathscr{W},\tau_\mathscr{V}$ have some components in common
\cite{GMMSRS}. \smallbreak

Recent work has focused on the problem of embedding classes of
$G$-structures on a parallelizable manifold inside an appropriate
parameter space. In this case, one can consider $G$-structures that
stabilize a global section $\xi$ of some tensor power of the tangent
bundle. Such structures are parametrized by a unique $G$-orbit
$\sO_\xi$. The \emph{intrinsic torsion varieties} (ITV's) of a
parallelizable manifold are the subsets of $\sO_\xi$ of structures
belonging to the same null-torsion class (see \cite{GMMSRS,GMSSTV}).
An analysis of ITV's of structures on the Iwasawa manifold and other
nilmanifolds has been carried out in \cite{AGSHGI,AGSAHG}. We
highlight some relevant examples which we combine with the techniques
of the present article.

The Iwasawa manifold $N$ is defined as the set of right cosets
$\Gamma \backslash G_H$, where $G_H$ is the complex Heisenberg group
and $\Gamma$ the natural lattice:
\[
G_H = \left\{ \left(\begin{array}{ccc}
1& z^1 & z^2 \\
0& 1 & z^3 \\
0& 0 & 1\end{array}\right):\,z^k\in\C\right\},\quad
\Gamma=\left\{\left(\begin{array}{ccc}
1& a^1 & a^2 \\
0& 1 & a^3 \\
0 &0 & 1
\end{array}\right):\,a^k\in \mathbb{Z}[i] \right\}\]
Nilmanifolds admit a natural parallelism determined by $G_H$-left
invariant vector fields. The complex $1$-forms $\xi_1=dz^1$,
$\xi_2=dz^2$ and $\xi_3=-dz^3+z^1dz^2$ are left invariant on $G_H$,
and can be exploited for defining the following left invariant real
$1$-forms:
\begin{equation}\label{baseiwa}\xi^1=e^1+\imath e^2\,\,\,\,\,\xi^2=e^3+\imath e^4\,\,\,\,\, \xi^3=e^5+\imath e^6\end{equation}

\noindent We will keep faithfully the notations of the previous
sections referring to this basis of the cotangent spaces of $N$.
Setting $e^i$ being orthogonal we define on $N$ a standard
Riemannian metric induced from a left-invariant tensor on $G_H$ (see
\cite{AGSHGI}). Consider then the $G$-structures obtained as
reductions of this fixed Riemannian structure by stabilizing a
$G_H$-left invariant $2$-form.\smallbreak

The following result has been proved as Theorem~1 in \cite{AGSHGI},
though we can now give it a moment map interpretation.

\begin{teor}(Abbena-Garbiero-Salamon) The set $I$ of invariant
  complex structures on $N$ is given by the disjoint union of the
  point $\omega_0$ and a $\CP^1$.  This is a $T$-invariant subset of
  $\sP^+$ and its image by $\mu_T$ is the union of a vertex and the
  edge $\mathcal{E}_{12}$ of $\Delta_{\sP^+}$.
\label{compiwa}
\end{teor}

\noindent A justification of the invariance of this set under the
action of a maximal torus, based upon the fact that $G_H$ is a
\emph{complex} Lie group, is given in the author's joint article
\cite{GMSSTV}, which also proves:

\begin{teor} The ITV of OPS's $P\in\sG$ on the Iwasawa manifold $N$
  characterized by an integrable 4-dimensional distribution (meaning
  $[\mathcal{H},\mathcal{H}] \subseteq\mathcal{H}$) is the complex
  submanifold
  \[ K= \Gr_2(\langle e^1,e^2,e^3,e^4\rangle)\cong \CP^1\times\CP^1\]
  of $\sG$ whose image $\mu_{T}(K)$ is the
  intersection $\Delta_{\sG}\cap \langle e^{12},e^{34}\rangle$.
\label{intvert}
\end{teor}

Let us denote by $K'$ the analogous subset of $\sG$ of OPS's
characterized by an integrable 2-dimensional distribution (meaning
$[\mathcal{V},\mathcal{V}]\subseteq\mathcal{V}$). From \cite{GMMSRS}
we know:

\begin{teor} The subset $K\cap K'$ of $\sG$ is the disjoint union of two
$2$-spheres consisting of the $J_0$-invariant elements of $K$. Its
image by the moment map is
\[\mu_T(K\cap K')=\{xe^{12}+ye^{34}:x+y=\pm1,\
|x|,|y|\le1\},\] and is formed of two line
segments.\label{doubint}\end{teor}

The OPS's in $K$ characterized by Theorem~\ref{intvert} are those for
which a Nihenjuis-type tensor $\Lambda^2\mH\to\mV$ is zero
\cite{NavAPS}, or equivalently those with vanishing $\mathscr{V}_4$
component in \eqref{VVVVVV}. The corresponding null-torsion class of
these \emph{foliations} is therefore $\mathscr{V}_{12356}$. A complete
analysis of the ITV's of OPS's on the Iwasawa manifold $N$ is given in
\cite{GMITCR}, as a result of which it turns out that
\begin{equation}\label{V12356}
\mathscr{V}_{12356}=\mathscr{V}_{15};
\end{equation}
indeed it is easy to see that elements of $K$ have OPS torsion in
$\mathscr{V}_1\oplus\mathscr{V}_5$.  Moreover, elements of $K\cap
K'$ describe those OPS's for which $\mH$ gives rise to a
\emph{totally geodesic foliation} and they are of class
$\mathscr{V}_5$. Another result from \cite{GMITCR} along these lines
regards a relevant subset of $K'$:

\begin{teor}\label{newexample4} The ITV of foliations of class
  $\mathscr{V}_{345}$ on the Iwasawa manifold $N$ is the disjoint
  union of two $\CP^2$'s in $\mathscr{G}$. The OPS's are
  characterized by the fact that their $2$-plane $\mV$ is
  $J_0$-invariant.
\end{teor}

\noindent The two $\CP^2$'s in question are precisely the subsets
previously denoted by $F^+$ and $F^-$.\medbreak

\begin{figure}[!h]

 \centering\includegraphics[width=0.9\textwidth]{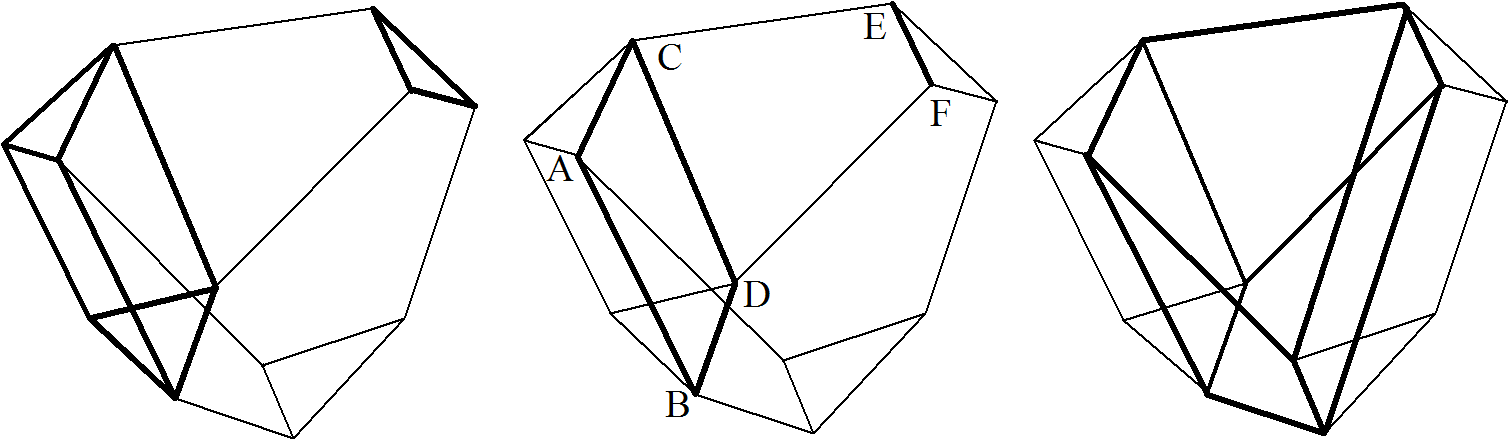}\\
 \caption{$\mu_T(\pi_1^{-1}(I)\cap\pi_2^{-1}(K))$}
 \label{mixed2}\end{figure}

Knowledge of the fibrations of coadjoint orbits described in
Section~1 enables us to detect ITV's of mixed structures inside any
10-dimensional orbit $SO(6)/U(1)\times U(2)$ (take for example
$\sF^+$ in \eqref{trianglefibr3}). Consider an OCS $J\in\sP^+$ and
an OPS $P\in\sG$ with intrinsic torsion in a known class.
Proposition \ref{nablajp} implies that the points of $\sF^+$ in the
intersection $\pi^{-1}_1(J)\cap\pi_2^{-1}(P)$ have pre-determined
intrinsic torsion. If we consider in this way the inverse images of
entire classes, their intersection determines a specific ITV of
mixed structures inside $\sF^+$. Theorem~\ref{compiwa} and Theorem~\ref{intvert} now yield:

\begin{corol}\label{cor1}
  The ITV inside $\sF^+$ consisting of MS's on $N$ of class
  $\mathscr{W}_{34}\cap\mathscr{V}_{15}$ is a disjoint union
  $(\CP^1\times \CP^1)\sqcup \CP^1$. It's image by the moment
  map is shown in the centre of Figure~\ref{mixed2}.
\end{corol}

\nit{Proof} By Proposition~\ref{square} $\pi_2^{-1}(K)$ is a real $6$-dimensional symplectic toric manifold, $\pi_2^{-1}(I)$ has two disjoint components, a $\CP^2$, and the symplectic toric manifold $\pi_1^{-1}(L)$ (recall Proposition~\ref{edge}). Proposition~\ref{refery} implies that the $C_i$-invariant subsets of $\pi_1^{-1}(L)$ and $\pi_2^{-1}(K)$ are real four-dimensional symplectic toric manifolds projecting on the faces of the corresponding moment polytopes, so the intersection of $\pi_1^{-1}(L)$ and $\pi_2^{-1}(K)$ is exactly the $C_2$-invariant subset projecting onto the common rectangular face. This means $\pi_1^{-1}(L)\cap\pi_2^{-1}(K)\cong \CP^2\times\CP^2$. A dimensional check shows that the subsets of $\pi_2^{-1}(K)$ and $\pi_2^{-1}(I)$ simultaneously invariant under the action of two $C_i$-s are all real two-dimensional symplectic manifolds, which project onto segments (edges of the polytopes), and are therefore toric manifolds symplectomorphic to $\CP^1$. The intersection of $\CP^2\subset\pi_2^{-1}(I)$ and $\pi_2^{-1}(K)$ is a $\CP^1$.\qed

\begin{corol}\label{cor2}
  The ITV inside $\sF^+$ consisting of MS's on $N$ of class
  $\mathscr{W}_{34}\cap\mathscr{V}_{5}$ is a disjoint union
  $\CP^1\sqcup  \CP^1\sqcup\CP^1$. The projection to $\Delta_{\sF^+}$
  consists of the segments $AB$, $CD$, $EF$ in Figure~\ref{mixed2}.
\end{corol}

\nit{Proof} This follows from Theorem~\ref{doubint}, the proof is analogous to the previous-one. Proposition~\ref{refery} implies that the trapezium faces of $\pi_2^{-1}(K)$ are projections of real four-dimensional $C_i$-invariant symplectic toric manifolds ($\CP^1$ bundles over $\CP^1$). \qed

\begin{figure}[!h]

\centering\includegraphics[width=0.9\textwidth]{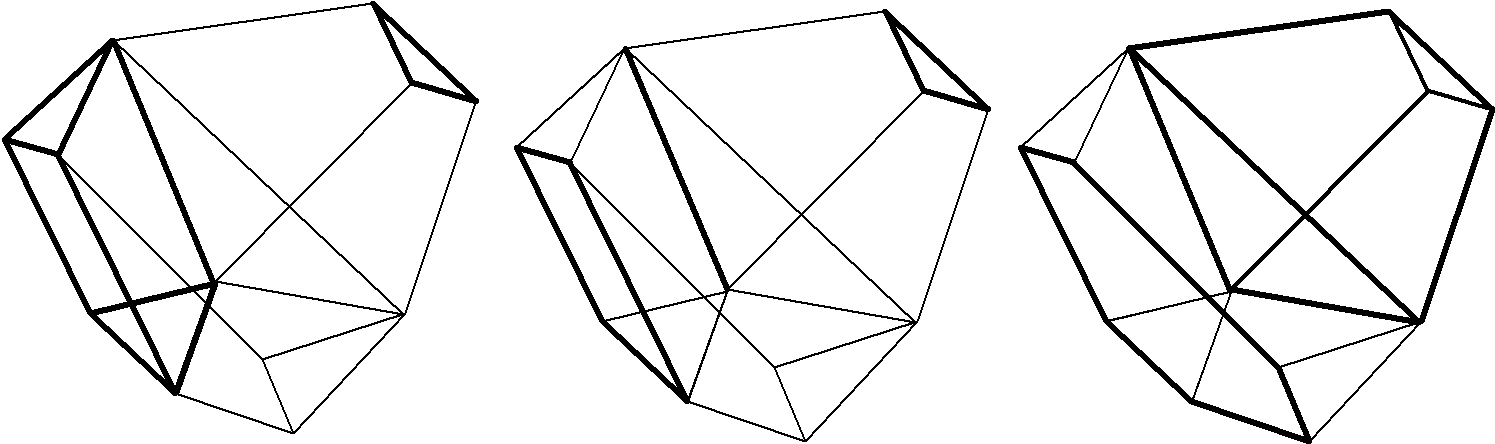}\\
\caption{$\mu_T(\pi_1^{-1}(I)\cap\pi_2^{-1}(F^+\cup F^-))$}
\label{newexamp}\end{figure}

\begin{corol}\label{cor3}
  The ITV inside $\sF^+$ consisting of MS's on $N$ of class
  $\mathscr{W}_{34}\cap\mathscr{V}_{345}$ is a disjoint union
  $\CP^1\sqcup\CP^2\sqcup(\CP^1\rtimes\CP^1)$. Its projection
  to $\Delta_{\sF^+}$ consists of a line segment, a triangle
  and a trapezium as shown in the centre of Figure~\ref{newexamp}.
\end{corol}\nit{Proof} This follows from Theorem~\ref{newexample4}. The case of $F^+$ is analogous to the previous cases. The trapezium face of $\mu_T(L)$ is the image of a $C_1$-invariant toric submanifold ($\CP^1$ bundle over $\CP^1$). Proposition~\ref{refery} implies that this is a four-dimensional toric submanifold of $F^-$.\qed

\begin{rem}The $C_i$-invariant (toric) submanifolds play a key role in the above results. Each corollary can be proved also in terms of $2$-forms and compatibility of $G$-structures. As an example we consider the intersection of $\pi_1^{-1}(L)$ and $\pi_2^{-1}(K)$ . Lemma~\ref{tec} implies that the set of MS's defined by an OCS in $L$ and an OPS in $K$ is a trivial $S^2$ bundle over $S^2$. It follows from \eqref{edges} that a generic $2$-form in $L$ has the
following expression$$\omega=a(e^{12}-e^{34})+b(e^{13}-e^{42})+c(e^{14}-e^{23})-e^{56}.$$
\noindent Without loss of generality, we can fix a unit $1$-form
$w=y_1 e^1+y_2 e^3+y_3 e^5$. Let $J$ be an OCS in $L$. Then
\[\begin{array}{c}
Jw=y_1(ae^2+be^3+ce^4)+y_2(-ae^4-be^1+ce^2)-y_3
e^6,\\[3pt] \mu_T(\omega+\alpha w\wedge Jw) =
\big(a+\alpha(y_1^2a+y_1y_2c),-a+\alpha(y_1y_2 c-y_2^2a),-1-\alpha
y_3^2\big).
\end{array}\]
Adding a $2$-form in $L$ to that of a consistently-oriented element
in $K$, we obtain a $2$-form with $e^{56}$ component equal to $-1$.
So $\mu_T(\omega+\alpha e\wedge Je)$ is the intersection of $\sF^+$
with the plane $z=-1$, etc... 
\end{rem}

\begin{rem} In the notation of \eqref{fine} and \eqref{MWV}, one
  can identify the subspaces of $\mathscr{M}$ containing the intrinsic torsion
  of the MS's described by Corollaries \ref{cor1},
  \ref{cor2} and \ref{cor3}. Namely:
\[\begin{array}{rcl}
\mathscr{W}_{34}\cap\mathscr{V}_5 &\cong&
\rr{\nu\la^{0,2}}\oplus\>\rr{\nu\la^{1,1}_0}\oplus\mH,\\[2pt]
\mathscr{W}_{34}\cap\mathscr{V}_{15} &\cong&
\rr{\nu\la^{0,2}}\oplus\>\rr{\nu\la^{1,1}_0}\oplus2\mH,\\[2pt]
\mathscr{W}_{34}\cap\mathscr{V}_{345} &\cong&
\rr{\nu\la^{0,2}}\oplus2\rr{\nu\la^{1,1}_0}\oplus2\mH.
\end{array}\]
Each class is therefore characterized by a relatively
small subset of the 16 irreducible
$U(1)\times U(2)$ components of $\mathscr{M}$.
\end{rem}

The last three corollaries show that the intersections of ITV's are
determined by the intersections of their moment polytopes. Actually,
it was the diagrams that led the author to formulate these results.
They should lead to a similar description of the intrinsic torsion
varieties for other structures on $N$ and on other nilmanifolds. The
techniques developed in this article should help establish the extent
to which these subsets of coadjoint orbits for $SO(6)$ are invariant
by tori. Except for the graphical aspects, the methods are not
restricted to the 6-dimensional case, but are perfectly general.

The Riemannian $G$-structures that we have considered can also be
obtained as reductions of \emph{spin} structures, since $Spin(6)$ is
isomorphic to $SU(4)$. Adopting a spinorial interpretation allows for
the possibility of enlarging the theory to include structures not
necessarily defined by a $2$-form. An important example in six
dimensions is that of $SU(3)$-structures whose intrinsic torsion
measures the extent to which a manifold fails to be Calabi--Yau. We
refer the reader to \cite{GMMSRS}.

\section*{Acknowledgements}
This article is based on part of the author's doctoral thesis
\cite{GMMSRS}, supervised by prof. Simon Salamon. The author thanks also the 
Geometriae Dedicata referee for the useful suggestions.

\small
\bibliographystyle{plain}
\bibliography{TMMRS}

\end{document}